\documentclass[letterpaper, 12pt]{article}

\title{Ideals of Adjacent Minors}
\author{Serkan Ho\c{s}ten  \\
{\small Department of Mathematics, 
San Francisco State University, San Francisco} \\
\and Seth Sullivant \thanks{Partially supported by the
National Science Foundation (DMS-0200729)} \\
{\small Department of Mathematics, University of California,
Berkeley} 
}
\date{}

\usepackage{epsfig}
\usepackage{latexsym,array,delarray,amssymb}

\newcommand{\comment}[1]{}

\newtheorem{thm}{Theorem}[section]
\newtheorem{lemma}[thm]{Lemma}

\newtheorem{cor}[thm]{Corollary}
\newtheorem{prop}[thm]{Proposition}

\newtheorem{question}[thm]{Question}

\newtheorem{Example}[thm]{Example}
\newtheorem{Remark}[thm]{Remark}
\newtheorem{Alg}[thm]{Algorithm}
\newtheorem{Defn}[thm]{Definition}

\newenvironment{example}{\begin{Example}\rm}
                {\mbox{}~\hfill~$\Box$~\end{Example}}

\newenvironment{defn}{\begin{Defn}\rm}{\end{Defn}}

\newenvironment{proof}{\begin{trivlist}\item {\it
        Proof.\,}}{\mbox{}~\hfill~$\Box$\end{trivlist}}

\newenvironment{proofof}[1]{\begin{trivlist}\item {\it
        Proof of {#1}.\,}}{\mbox{}\hfill$\Box$\end{trivlist}}

        {\begin{list}
                {\noindent\makebox[0mm][r]{(\roman{enumi})}}
                {\usecounter{enumi} \topsep=1.5mm \itemsep=0mm}
                
        }
        {\end{list}}

        {\begin{list}
                {\noindent\makebox[0mm][r]{\arabic{enumi}.}}
                {\usecounter{enumi} \topsep=1.5mm \itemsep=-1mm}
        }
        {\end{list}}

\def\<{\langle}
\def\>{\rangle}
\def\0{{\mathbf 0}}
\def\1{{\mathbf 1}}
\def\AA{{\mathbb A}}

\def\CC{{\mathbb C}}

\def\QQ{{\mathbb Q}}

\def\ZZ{{\mathbb Z}}

\def\xx{{\mathbf x}}

\def\IN{\mathsf{in}}

\def\gl{{G_{\!}L}}

\begin{document}

\maketitle

\begin{abstract}
We give a description of the minimal primes of the ideal generated
by the $2 \times 2$ adjacent minors of a generic matrix. We also
compute the complete prime decomposition of the ideal of adjacent
$m \times m$ minors of an $m \times n$ generic matrix when the
characteristic of the ground field is zero.  A key
intermediate result is the proof that the ideals which appear as
minimal primes are, in fact, prime ideals.  This introduces a
large new class of mixed determinantal ideals that are prime.
\end{abstract}

\section{Introduction}
Let $X_{mn}$ be an $m \times n$ matrix of indeterminates $x_{ij}$
which generate the polynomial ring $K[x_{ij}]$ where $K$ is a field.
The ideal generated by all $k \times k$ minors of $X_{mn}$
has been studied from many
different points of view; for a comprehensive exposition see \cite{BV}
and \cite[Chapter 7]{BH}. For
example, these ideals are prime ideals that are also Cohen-Macaulay
\cite{HE}, and they are Gorenstein when $m=n$ \cite{S}.
Similar determinantal ideals where one mixes minors
of different sizes have been also studied. For instance, in the
context of invariant theory and algebras with straightening
laws one looks at the ideal of minors generated by a
coideal in a particular poset of all minors \cite{DEP}.
There are also many variations such as {\em ladder determinantal
ideals} \cite{Conca}, and {\em mixed ladder determinantal ideals}
\cite{GM} where the ideals of (mixed) minors in a ladder-shape
region in $X_{mn}$ are studied. In both cases these ideals are prime
and Cohen-Macaulay, and criteria for when they are Gorenstein are
characterized.

A $k \times k$ {\em adjacent minor} of $X_{mn}$ is the determinant
of a submatrix with row indices
$r_1, \ldots, r_k$ and column indices $ c_1, \ldots, c_k$
where these indices are consecutive integers.
We let $I_{mn}(k)$ be the ideal
generated by all of the $k \times k$ adjacent
minors of $X_{mn}$. As opposed to the ideal of {\em all}
$k \times k$ minors, the ideal $I_{mn}(k)$ is far from being a
prime ideal. This ideal first appeared in \cite{DES}
for the case $k = 2$  where primary decompositions
of $I_{2n}(2)$ and $I_{44}(2)$ were given. The motivation
for studying $I_{mn}(2)$ comes from the rapidly
growing field of {\em algebraic statistics} \cite{PRW},
\cite[Chapter 8]{St}: a primary decomposition of $I_{mn}(2)$ helps
to measure the connectedness of the set of $m~\times~n$
{\em contingency tables} with the same row and column sums
via the moves corresponding to the $2 \times 2$ adjacent minors
\cite{DES}.

The goal of this paper is to study the minimal primes of $I_{mn}(k)$.
A motivation  is related to
algebraic statistics and focuses on the case when $k=2$ in
Section 2, and on the case of adjacent minors of higher-dimensional
matrices in Section 5. We give in Section 2 a combinatorial
description of the minimal primes of $I_{mn}(2)$.
This ideal is a very special instance of a {\em lattice basis ideal},
and minimal primes of lattice basis ideals have been characterized
\cite{HS}. However, in the case we treat here we get a more
transparent characterization.

In Section 3 we analyze the case when $k = m$, i.e. the maximal
adjacent minors of an $m \times n$ matrix where $m \leq n$. In
this case, $I_{mn}(m)$ is a complete intersection that is also
radical.  We present a combinatorial description of the minimal
primes and give a recurrence relation for the number of these primes.
These prime ideals are a very general type of mixed determinantal
ideals that, to our knowledge, have never before been studied. All
the usual questions can be asked about them, however, even the
fact that they are prime seems to be a challenging  result.  Section 4
is the technical heart of the paper: it is devoted to the proof
that these mixed determinantal ideals are, in fact, prime.  A
string of arguments that culminates in Theorem \ref{thm:prime}
proves this result when $\mathrm{char}(K) = 0$. 
In arbitrary characteristic we also show that they
are prime in special cases including the case 
when $m \leq 3$.  On the way
to proving these results we show that the minors that generate
these mixed determinantal ideals form a squarefree Gr\"obner
basis when the characteristic is arbitrary.

Section 5 is a look into the future with a view towards
applications in algebraic statistics. We introduce the notion of
adjacent minors of a generic $m_1 \times m_2 \times \cdots \times
m_d$ matrix. These come from the study of discrete random
variables $X_1, \ldots, X_d$ where each $X_i$ takes values in
$\{1, \ldots, m_i\}$. A particular family of statistical models
that describe the joint probability distributions of these random
variables (the so-called no $d$-way interaction models \cite{F}) 
gives rise
to a toric variety whose set of defining equations may be
extremely large and complicated \cite{AT, St}. 
However, the {\em positive} probability
distributions are described precisely by the simple
multidimensional adjacent minors we will introduce.  The story of
the minimal primes of these ideals is far from complete, but in
Theorem \ref{thm:22n} we will describe them in the case $m_1 = m_2
= \cdots = m_{d-1} = 2$.

\section{$2 \times 2$ Adjacent Minors}

From the general characterization of minimal primes of lattice
basis ideals \cite{HS} it follows  that every minimal prime $P$ of
$I_{mn}(2)$ is of the form
\begin{equation} \label{primeform}
P = \< x_{ij} \, : \,\, x_{ij} \in S \> \,\, +  \,\,\,
J : (\prod_{x_{ij} \notin S} x_{ij})^\infty
\end{equation}
where $S$ is a subset of the variables in the ring $K[x_{ij}]$ and
$J$ is the ideal generated by the $2 \times 2$ adjacent minors
in the ring $K[x_{ij}: \,\, x_{ij} \notin S]$. In other words,
$P$ is uniquely determined by the variables it contains. We will
denote this set of variables by $S_P$, and the variables not in
$S_P$ by $N_P$. In the rest of this section we will
give a characterization of the sets $S_P$ and $N_P$ that give
rise to the minimal primes of  $I_{mn}(2)$. In order to describe
these minimal primes we need a few definitions.

Let $S$ be a subset of variables of $K[x_{ij}]$.  We say that
two variables $x_{ij}$ and $x_{st}$ are \emph{adjacent} if
$s = i +\epsilon_1$ and $t = j + \epsilon_2$ where
$\epsilon_1, \epsilon_2 \in \{-1,0,1\}$.  The set $S$ is
\emph{connected} if for every pair of variables
$\{x_{ij}, x_{st} \} \subset S$ there is a sequence of variables in $S$
starting with $x_{ij}$ and  ending with $x_{st}$, and such that
each variable in the sequence is adjacent to the variable preceding and
following it.  A subset $T$ of $S$ is called {\em maximally connected} if
there is no larger connected subset of $S$ containing $T$.
A set of variables $S$ is a {\em rectangle}
$X[i,j;s,t]$ if it is equal to the set of all the
variables in the submatrix
$$
\left(
\begin{array}{ccc}
x_{ij} & \cdots & x_{it} \\
\vdots & \ddots & \vdots \\
x_{sj} & \cdots & x_{st}
\end{array}
\right).
$$
The {\em boundary edges} of $X[i,j;s,t]$ are the four rectangles
$X[i-1,j;i-1,t]$, $X[s+1,j;s+1,t]$,
$X[i,j-1;s,j-1]$, and $X[i,t+1;s,t+1]$.
The {\em boundary} of $X[i,j;s,t]$ is the union of the
the four boundary edges together with the ``corner'' variables
$x_{i-1,j-1}$, $x_{s+1,j-1}$, $x_{i-1,t+1}$, and
$x_{s+1,t+1}$. When we speak of boundary edges
and the boundary of a rectangle we always
mean only those parts that are defined, since some boundary
edges or corner variables might not exist because they
are outside of the matrix $X_{mn}$.

\begin{figure} \label{fig1}
\centerline{
\epsfig{file=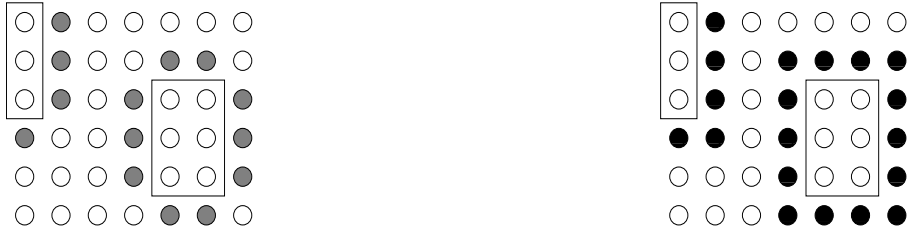, height={3cm}, width={12cm}}}
\caption{Two rectangles with their boundary edges and boundaries}
\end{figure}

\begin{example}
Let $m = 6$ and $ n = 7$. In the matrix $X_{67}$, the two rectangles
$X[1,1;3,1]$ and $X[3,5;5,6]$ together with their
boundary edges and boundaries can be viewed in Figure 1.
The first rectangle has only two boundary edges since the
other two are not defined.
\end{example}

\begin{defn} \label{primepartition}
We will call a partition $(S,N)$ of the variables in $X_{mn}$
a {\em prime partition}  if $S$ and $N$ satisfy the following
properties:
\begin{enumerate}
\item $N$ contains the variables $x_{11}$, $x_{1n}$, $x_{m1}$
and $x_{mn}$,
\item when $N$ is written as the disjoint union of  its maximally
connected subsets $N = \bigcup_k T_k$, then each $T_k$ is a rectangle,
\item each boundary edge of a maximal rectangle $T_k$ in $N$ has a
nontrivial intersection with the boundary of another maximal rectangle
$T_{\ell}$,
\item the boundary edges of two maximal rectangles of
width (height) one in the same column (row) do not intersect, and
\item $S$ is the union of the boundaries of the maximal rectangles $T_k$.
\end{enumerate}
\end{defn}

\begin{thm} \label{primethm}
The prime ideal $P$ is a minimal prime of $I_{mn}(2)$ if and
only if $(S_P, N_P)$ is a prime partition.
\end{thm}

The rest of the section is devoted to the proof of Theorem
\ref{primethm}. We remark that this theorem does indeed cover
the characterizations of minimal primes of $I_{mn}(2)$ in the
known cases, in particular, that of $I_{2n}(2)$ in \cite{DES} and
of $I_{3n}(2)$ in \cite{HS}. Before starting the proof we give an
example to illustrate the definition above and the content of the theorem.

\begin{figure} \label{fig2}
\centerline{
\epsfig{file=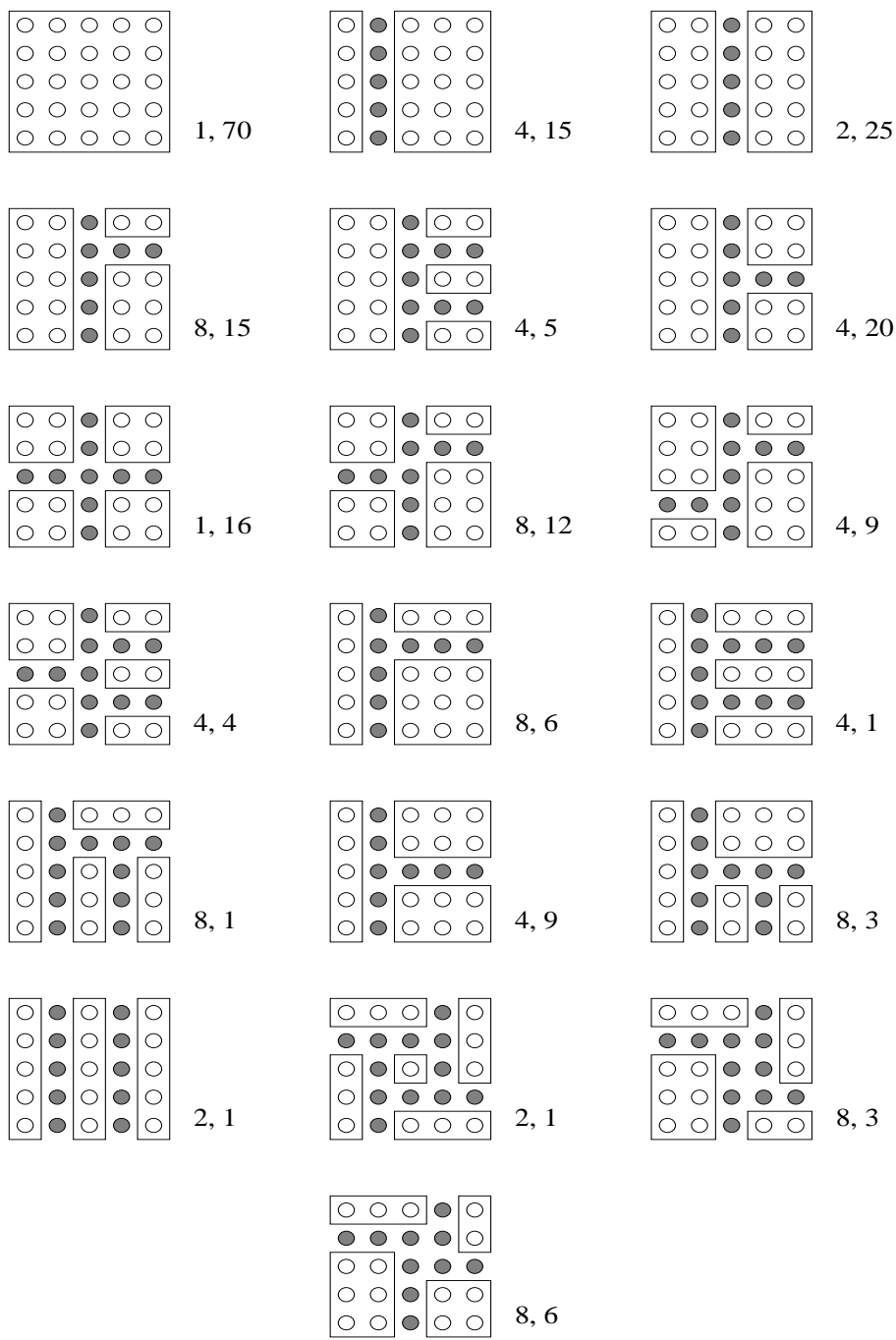, height={18cm}, width={12cm}}}
\caption{The minimal primes of $I_{55}(2)$}
\end{figure}

\begin{example}
Figure \ref{fig2} displays all the minimal primes of $I_{55}(2)$.
This is the smallest example where all five conditions in the
Definition \ref{primepartition} are needed.  In this case there
are 92 minimal primes that can be grouped into 19 equivalence
classes modulo symmetries. We show one member from each
equivalence class. The boxes in Figure \ref{fig2} are the maximal
rectangles in the $N_P$ of the corresponding prime partition, and
the solid buttons correspond to the variables in $S_P$. The first
number following each diagram is the size of the equivalence class
and the second is the degree of the corresponding prime ideal.
\end{example}

Now we begin the proof of Theorem \ref{primethm} with a sequence
of lemmas. The first one concerns the first property in
Definition \ref{primepartition} and is taken from Lemma 3.3 in \cite{HS}.

\begin{lemma} \label{lem:corners}
The corner variables $x_{11}$, $x_{1n}$, $x_{m1}$, and $x_{mn}$ do
not belong to $S_P$ for any minimal prime $P$ of $I_{mn}(2)$.
\end{lemma}

\begin{lemma}\label{lem:rectangle}
If $P$ is a minimal prime of $I_{mn}(2)$,  then every
maximally connected subset of $N_P$ is a rectangle.
\end{lemma}
\begin{proof}
Let $T$ be a maximally connected subset of $N_P$ and suppose
that the adjacent variables $x_{ij}$ and $x_{i+1,j+1}$ are in $T$.
Since these two variables are not in $P$, the only way
the adjacent minor $x_{ij}x_{i+1,j+1}-x_{i,j+1}x_{i+1,j}$
could be in $P$ is if the variables $x_{i+1,j}$ and $x_{i,j+1}$ also
belong to $N_P$. Since $T$ is maximally connected these
two variables are also in $T$.  Similarly, if $x_{i+1,j}$ and
$x_{i,j+1}$ belong to $T$ then $x_{ij}$ and $x_{i+1,j+1}$ are 
also in $T$.  This implies that any maximally connected
subset of $N_P$ is a rectangle.
\end{proof}

The general description of the minimal primes in (\ref{primeform})
together with Lemma \ref{lem:rectangle} imply that if $P$ is a
minimal prime of $I_{mn}(2)$, and $N_P$, the set of variables not
in $P$, is written as the disjoint union of its maximally
connected rectangles, say $N_P = \bigcup_k T_k$, then
$$ P = \< x_{ij}: \,\, x_{ij} \in S_P \> \, + \,
\<x_{ij}x_{st} - x_{it}x_{sj}: \mbox{ all of } x_{ij}, x_{st},
x_{it}, x_{sj} \mbox{ are in the same } T_k \>.$$

\begin{lemma}\label{lem:boundary}
Let $P$ be a minimal prime of $I_{mn}(2)$ and let $T$ be a maximally
connected rectangle of $N_P$. Then the boundary of $T$
is a subset of $S_P$. Moreover, for each boundary edge $E$ of $T$
there is another maximal rectangle $T' \subset N_P$ whose
boundary has a nonempty intersection with $E$.
\end{lemma}

\begin{proof}
The boundary of $T$ is a subset of $S_P$ since $T$ is maximally
connected.  To prove the second statement,  suppose
that there were a maximal rectangle $T$ with a boundary edge $E$ that
does not intersect the boundary of any other maximal rectangle.
Consider the prime ideal $P'$ where $S_{P'}= S_P \setminus E$,
and $N_{P'} = N_{P} \cup E$. The assumption on the edge $E$ implies
that $T' = T \cup E$ is a maximally connected rectangle of $N_{P'}$.
The new prime ideal $P'$ still contains all the adjacent minors.
The only new $2 \times 2$ minors that
appear in the ideal $P'$ involve variables from $E$, and
these are already contained in $P$. This implies that $P'$ is a
prime ideal contained in $P$, contradicting
the minimality of $P$.
\end{proof}

\begin{lemma}\label{lem:vector}
Let $P$ be a minimal prime of $I_{mn}(2)$ and let $T = X[i,j;i,s]$ be a
maximally connected rectangle in $N_P$ of height one. Then there
is no maximally connected rectangle of height one in $N_P$ of the
form $T' = X[i,s+2;i,t]$. A similar statement holds for
vertical rectangles of width one.
\end{lemma}
\begin{proof}
By Lemma \ref{lem:boundary} the rectangles
$X[i-1,j;i-1,t]$ and $X[i+1,j;i+1,t]$, and the variable $x_{i,s+1}$
are in $S_P$. Since the variables of $T$
and $T'$ do not appear in any generator of $P$, the prime ideal
$P'$ given by the set of variable $S_P \setminus x_{i,s+1}$ is a
strictly smaller prime ideal which contains $I_{mn}(2)$,
contradicting the minimality of $P$.
\end{proof}

\begin{lemma}\label{lem:boundary2}
If $P$ is a minimal prime of $I_{mn}(2)$, then every variable in
$S_P$ belongs to the boundary of some maximal rectangle in $N_P$.
\end{lemma}
\begin{proof}
Suppose that $S_P$ contains a variable $x_{ij}$ that is not in the
boundary of any maximal rectangle in $N_P$.  This implies that $x_{ij}$
is adjacent only to variables in $S_P$.  Let $U \subset S_P$
be the set of variables that are adjacent to $x_{ij}$.
The $2 \times 2$ adjacent minors contained in the ideal generated
by the variables in $U \cup x_{ij}$ are the same as those
contained in the ideal generated by the variables in $U$ alone.
Hence by omitting the variable $x_{ij}$ from $P$ we can construct
a prime ideal that contains $I_{mn}(2)$, but strictly contained in
$P$.  This is a contradiction to the minimality of $P$.
\end{proof}

With the help of the five lemmas we have presented we are ready to prove
the main theorem of this section.

\begin{proofof}{Theorem \ref{primethm}}
If $P$ is a minimal prime of $I_{mn}(2)$, the partition
$(S_P, N_P)$ satisfies all the five properties to be
a prime partition because of the five lemmas, Lemma \ref{lem:corners}
through Lemma \ref{lem:boundary2}, above. Hence we just need
to prove the converse. Suppose $(S, N)$ is a prime partition, and
we assume $N = \cup_k T_k$ is the partition of $N$ into its
maximally connected rectangles. We will show that the prime ideal
$$ P = \< x_{ij}: \,\, x_{ij} \in S \> \, + \,
\<x_{ij}x_{st} - x_{it}x_{sj}: \mbox{ all of } x_{ij}, x_{st},
x_{it}, x_{sj} \mbox{ are in the same } T_k \>$$ is a minimal prime of
$I_{mn}(2)$.  Since all the $T_k$ are rectangles, it is easy to
see that $P$ contains $I_{mn}(2)$. Suppose that there were a
minimal prime $P'$ over $I_{mn}(2)$ strictly contained in $P$.
This means that $(S_{P'}, N_{P'})$ is a prime partition, and
$S_{P'}$ is a proper subset of $S_P = S$. We consider a variable
$x_{ij}$ in $S_{P} \setminus S_{P'}$.  By Lemma
\ref{lem:boundary2}, $x_{ij}$ lies on the boundary of some maximal
rectangle $T$ of $N_P = N$. The variable $x_{ij}$ either lies  on
a boundary edge $E$ of $T$, or is a corner variable on the
boundary of $T$. In the first case, since $(S_{P'}, N_{P'})$ is a
prime partition, $E \subset S_P \setminus S_{P'}$, and therefore
$E$ is a subset of $N_{P'}$.  Moreoever $E$ intersects the
boundary of at least one other rectangle $T'$ of $N_P$.  This
means $T \cup E \cup T'$ is a connected subset of $N_{P'}$, and
this union must be contained in a maximally connected rectangle
$T''$ of $N_{P'}$. If $x_{ij}$ is a corner variable of the
boundary of $T$, then the two boundary edges $E$ and $E'$ of $T$
that are adjacent to $x_{ij}$ must be a part of $N_{P'}$. Now by
repeating the above argument we are guaranteed to have another
rectangle $T'$ of $N_{P}$ where $T \cup E \cup T'$ is contained in
a maximally connected rectangle $T''$ of $N_{P'}$. By the fourth
property of Definition \ref{primepartition}, $T$ and $T'$ could
not be both height (width) one rectangles in the same row (column)
of $X_{mn}$. Hence there are  variables $x_{st} \in T$ and $x_{pq}
\in T'$ where $s \neq p$ and $t \neq q$. Since these variables are
in the same maximally connected rectangle $T''$ of $N_{P'}$, the
$2 \times 2$ minor $x_{st}x_{pq}-x_{sq}x_{pt}$ is in $P'$. On the
other hand, the set of variables appearing in this minor is not
contained in any maximally connect rectangle of $N_P$ and so it
does not belong to $P$.  This contradicts the assumption that $P'
\subset P$.
\end{proofof}

There are many open questions left to answer about $I_{mn}(2)$.  A
combinatorial description of the embedded primes remains elusive.
Moreover, there are many interesting open questions regarding the
minimal primes.  For example, how many are there, which minimal
primes have the largest dimension, and what is the degree of the
radical $ \mathrm{rad}(I_{mn}(2))$?

\section{Maximal Adjacent Minors}

In this section we will describe the complete primary
decomposition of the ideals $I_{mn}(m)$ for $m \leq n$ over a
field $K$ of characteristic zero, and for
$m \leq 3$ in arbitrary characteristic.  With no restrictions on
the characteristic of the field
our description presents $I_{mn}(m)$ as the
irredundant intersection of radical ideals.

\begin{prop} \label{complete-inter}
The ideal $I_{mn}(m)$ is a radical ideal that is a complete
intersection.  Its codimension is $n-m+1$ and it has degree
$m^{n-m+1}$.
\end{prop}
\begin{proof} With respect to the lexicographic term order
where $x_{11} \succ x_{12} \succ \cdots \succ x_{1n} \succ
x_{21} \succ \cdots \succ x_{mn}$, the set of $m \times m$
adjacent minors of $X_{mn}$ is a Gr\"obner basis of $I_{mn}(m)$.
This follows from the fact that the initial terms of these
minors are pairwise relatively prime. The initial
ideal is a radical ideal that is a complete intersection,
and hence so is $I_{mn}(m)$.
Since there are $n-m+1$ maximal adjacent minors, the
codimension of $I_{mn}(m)$ is $n-m+1$ and its degree is $m^{n-m+1}$.
\end{proof}

Below we will give a description of the minimal primes of
$I_{mn}(m)$.  In this section we will show that $I_{mn}(m)$ is the
irredundant intersection of these radical ideals.  The proof that they
are prime in characteristic zero and when $m \leq 3$ for
arbitrary characteristic occupies Section 4.

\subsection*{Description of the minimal primes}

In order to make the narrative cleaner we will assume that the
matrix $X_{mn}$ has two {\em phantom} columns: a column indexed by
$0$ and another by $n+1$. (The role of the phantom columns is
only to make the description of the minimal primes simpler.) We will
denote by $[i,j]$ with $0 \leq i \leq j \leq n+1 $ the interval of
column indices $\{i,i+1, \ldots, j-1, j\}$ of $X_{mn}$, and
$X[i,j]$ will denote the submatrix consisting of the corresponding
columns of $X_{mn}$.

\begin{defn}
Let $\Gamma = \{[a_1,b_1],[a_2,b_2], \ldots, [a_k,b_k]\}$ be
a sequence of $k$ intervals. The sequence $\Gamma$ is called a
{\em prime sequence} if it satisfies the following properties:

\begin{enumerate}
\item $\bigcup [a_i,b_i] = [0,n+1]$,
\item $a_i < a_{i+1}$, $b_i < b_{i+1}$ for all $i$,
\item $b_i - a_i > m$ for all $i$, and
\item $0 \leq b_i - a_{i+1} < m-1$ for all $i$.
\end{enumerate}
\end{defn}

\noindent The definition says that each interval of $\Gamma$ is a
block of more than $m$ columns and all together they cover all the
columns of $X_{mn}$ (including the two phantom columns). Moreover the
consecutive intervals in the sequence have a nonempty overlap of
width less than $m$.   Given a prime sequence $\Gamma$ we let
$P_\Gamma$ be the ideal in $K[x_{ij}]$ defined by
\begin{enumerate}
\item all $m \times m$ minors of  $X[a_i,b_i]$
for each $[a_i,b_i] \in \Gamma$, and
\item all (maximal) $(b_{i}-a_{i+1} +1) \times (b_{i}-a_{i+1} +1)$ minors
of  $X[a_{i+1}, b_i]$ for $1 \leq i \leq k-1$.
\end{enumerate}
In other words,
$P_\Gamma$ is generated by the $m \times m$ minors of the
submatrices whose columns are indexed by the intervals in $\Gamma$, and
the maximal minors of the submatrices whose columns are indexed
by the overlap of consecutive intervals. An example will do the
best job to illustrate this construction.

\begin{figure} \label{fig3}
\centerline{
\epsfig{file=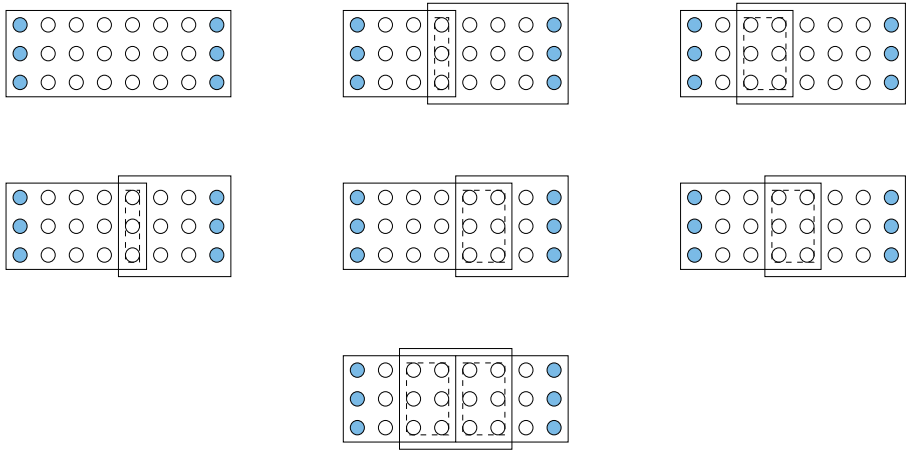, height={6cm}, width={12cm}}}
\caption{The minimal primes of $I_{36}(3)$}
\end{figure}

\begin{example}
We display the minimal primes $P_\Gamma$ of $I_{36}(3)$.  There
are seven primes corresponding to the seven prime sequences:
$$\begin{array}{rcl}
\Gamma_1 & = & \{[0,7]\} \\
\Gamma_2 & = & \{[0,3],[3,7]\} \\
\Gamma_3 & = & \{[0,3],[2,7]\} \\
\Gamma_4 & = & \{[0,4],[4,7]\} \\
\Gamma_5 & = & \{[0,4],[3,7]\} \\
\Gamma_6 & = & \{[0,5],[4,7]\} \\
\Gamma_7 & = & \{[0,3],[2,6],[4,7]\}
\end{array}
$$
Figure \ref{fig3} illustrates these minimal primes. The rectangles
with the solid borders describe the intervals in the corresponding
prime sequence. All $3 \times 3$ minors of each rectangle are
included in the corresponding minimal prime. We also indicate the
overlaps by rectangles with dashed borders; all the maximal minors
in these submatrices also need to be included in the corresponding
minimal prime.
\end{example}

\subsection*{The Main Theorem}
We now present the proof that the ideals $P_\Gamma$ describe the prime
decomposition of $I_{mn}(m)$ in characteristic zero and when $m
\leq 3$.  The following lemma will be needed for the proof of
Theorem \ref{thm:main}.
\begin{lemma} \label{lem:algo}
The variety $\mathcal{V}(I_{mn}(m))$ is contained in  
$\bigcup \mathcal{V}(P_{\Gamma})$ where the union is taken
over all prime sequences of $[0,n+1]$.
\end{lemma}
\begin{proof}
We will show that for each matrix 
$\mathbf{X} \in \mathcal{V}(I_{mn}(m))$ there is
a prime sequence $\Gamma$ such that $\mathbf{X}
\in \mathcal{V}(P_{\Gamma})$.
We describe an algorithm that constructs this prime sequence $\Gamma$.
For this,
let $\mathcal{I}(\mathbf{X}) := \{[c_1,d_1], \ldots, [c_t, d_t]\}$
be the set of all intervals of width less than $m$ in $[1,n]$ such that
$\mathbf{X}[c_i,d_i]$ has rank $d_i-c_i$, and $[c_i, d_i] \not \subset
[c_j, d_j]$ for $i \neq j$. We assume that $c_1 < \cdots < c_t$.
We define a prime sequence $\Gamma$ as follows:

\begin{enumerate}
\item Set $i = a_1 = b_0 = 0$ and $\Gamma = \emptyset$.
\item While $b_i \neq n+1$ do
\begin{itemize}
\item[(a)] i := i + 1.
\item[(b)] Let $[c_{j_i}, d_{j_i}] \in \mathcal{I}(\mathbf{X})$ be the
first interval in $[a_i, n+1]$ with $c_{j_i} > a_i + 1$. If there
is no such interval set $b_i = n + 1$.
\item[(c)] If $d_{j_i} \leq a_i + m$ set $b_i = a_i +m$,
unless $a_i + m \geq n$ in which case set $b_i = n+1$. Otherwise
set $b_i = d_{j_i}$.
\item[(d)] $\Gamma := \Gamma \cup \{[a_i, b_i]\}$.
\item[(e)] If $b_i \neq n +1$, let $[p_{j_i}, q_{j_i}] \in
\mathcal{I}(\mathbf{X})$ be the last interval in $[a_i,b_i]$. Set
$a_{i+1} = p_{j_i}$.
\end{itemize}
\item If the last interval in $\Gamma$ has width less than $m+1$ replace
it with $[n+1-m, n+1]$.
\end{enumerate}
Step $2\mathrm{(c)}$ together with step $3$ guarantees that the
intervals in $\Gamma$ have width at least $m+1$. Moreover, step
$2\mathrm{(e)}$ implies that consecutive intervals have a nonempty
overlap of width less than $m$. These show that $\Gamma$ is a prime
sequence.

Next we show that $\mathbf{X}$ is in $\mathcal{V}(P_\Gamma)$.
By the above construction of $\Gamma$ the overlap
$[a_{i+1}, b_i]$ of two
consecutive intervals contains one of the elements $[c_i,d_i]$ of
$\mathcal{I}(\mathbf{X})$. Since $\mathbf{X}[c_i,d_i]$ is
rank-deficient (it has rank $d_i-c_i$ instead of $d_i-c_i+1$),
so is $\mathbf{X}[a_{i+1},b_i]$, and the corresponding
$(b_i-a_{i+1}+1) \times (b_i-a_{i+1}+1)$ minors vanish on $\mathbf{X}$.

We need to show that the rank of $\mathbf{X}[a_i,b_i]$ for each
$[a_i,b_i] \in \Gamma$ is at most $m-1$. For this we analyze a few
different cases. First suppose that the width of $[a_i,b_i]$ is
bigger than $m+1$. The above algorithm implies that there
are either zero, one, two, or three intervals from
$\mathcal{I}(\mathbf{X})$ that are in $[a_i, b_i]$. When there are
no such intervals then $\Gamma = \{[0,n+1]\}$, and the matrix
$\mathbf{X}$ does not have any rank-deficient submatrices
consisting of less than $m$ adjacent columns. So
$\mathbf{X}[1,m-1]$ has full rank and these columns generate a
subspace $V$ with $\dim(V) = m-1$. But since the span of
$\mathbf{X}[1,m]$ is also $V$ and $\mathbf{X}[2,m]$ has rank
$m-1$, the span of $\mathbf{X}[2,m+1]$ and hence the span of
$\mathbf{X}[1,m+1]$ is $V$. Now by induction it is easy to see
that the span of $\mathbf{X}$ is the $(m-1)$-dimensional space
$V$, and therefore all $m \times m$ minors vanish on $\mathbf{X}$.
If there is one interval from $\mathcal{I}(\mathbf{X})$ inside
$[a_i, b_i]$, then either $[a_i,b_i] = [a_i, n+1]$ and the only
minimal rank-deficient interval is of the form $[a_i,c]$ with $c <
n$, or $[a_i, b_i] = [0,b_i]$ and the only minimal rank-deficient
interval is of the form $[c,b_i]$ with $c>1$. In the first case,
the submatrix $\mathbf{X}[a_i+1,n]$ has at least $m$ columns, and
this matrix does not have any rank-deficient submatrices
consisting of less than $m$ adjacent columns. By the same
argument above we conclude that the span of $\mathbf{X}[a_i+1,n]$
is an $(m-1)$-dimensional subspace $V$. But since
$\mathbf{X}[a_i,c]$ is minimally rank-deficient we conclude that
the span of $\mathbf{X}[a_i,n]$ is $V$, and therefore all $m
\times m$ minors corresponding to this interval vanish on
$\mathbf{X}$.  A symmetric argument applies when $[a_i, b_i] =
[0,b_i]$.

In the case where there are two intervals from
$\mathcal{I}(\mathbf{X})$, the two minimally rank-deficient
intervals are of the form $[a_i, c]$ and $[d, b_i]$ where $c <
b_i$ and $d > a_i$ \emph{or} of the form $[a_i,c]$ and $[a_i+1,
d]$ which forces the interval $[a_i, b_i] = [a_i,n+1]$.  This
means that, in the first case, $\mathbf{X}[a_i+1,b_i-1]$ has at
least $m$ columns and does not have any rank-deficient submatrices
consisting of less than $m$ adjacent columns. Similar
considerations as above show that $\mathbf{X}[a_i+1,b_i-1]$ has
rank $m-1$.  Since $\mathbf{X}[a_i,c]$ is minimally rank deficient
and $\mathbf{X}[a_i+1,a_i+m]$ is not rank deficient we see that
the column of $\mathbf{X}$ indexed by $a_i$ is in the span of the
columns of $\mathbf{X}[a_i+1,a_i+m]$ and so $\mathbf{X}[a_i,b_i]$
has rank $m-1$.  In the second case, the usual argument implies
that $\mathbf{X}[a_i +1, n]$ has rank $m-1$.  But since
$\mathbf{X}[a_i,c]$ is minimally rank deficient and
$\mathbf{X}[a_i+1,c]$ is not rank deficient we see that the column
of $\mathbf{X}$ indexed by $a_i$ is in the span of the columns of
$\mathbf{X}[a_i+1,n]$ and so $\mathbf{X}[a_i,n]$ has rank
$m-1$.  Finally, we consider the case where there are three
intervals from $\mathcal{I}(\mathbf{X})$ in $[a_i, b_i]$.  By
construction these are necessarily of the form $[a_i,c], [a_i +1,
d]$, and $[e, b_i]$.  But then the combination of the two
arguments for the cases with two minimally rank deficient
intervals shows that $\mathbf{X}[a_i,b_i]$ has rank $m-1$.

The case where the width of $[a_i,b_i]$ is exactly $m+1$ requires
a slightly different argument. If $[a_i,b_i]=[0,m]$ or $[n+1-m, n+1]$
there is nothing to show since there is only one $m \times m$ minor
that needs to be considered and it is necessarily an adjacent minor.
If we are not in these two trivial cases, the construction of $\Gamma$
implies that there are at least two intervals from
$\mathcal{I}(\mathbf{X})$ contained in $[a_i, b_i]$. Let $[c,d]$ be
the first such interval and $[e,f]$ the last such interval. Observe
that we have $c = a_i$. Now if these two intervals do not
overlap then any $m \times m$ submatrix of $\mathbf{X}[a_i,b_i]$
will contain one of these rank-deficient intervals and hence its
rank will be at most $m-1$. If there is an overlap we have
$a_i < e \leq d < f \leq b_i$. The rank of the submatrix
$\mathbf{X}[a_i,d]$ is $d-a_i$, and the rank of $\mathbf{X}[e,f]$ is
$f-e$. Moreover, since these intervals are minimally rank-deficient 
the rank of $\mathbf{X}[e,d]$ is $d-e+1$. But then the rank
of $\mathbf{X}[a_i,b_i]$ is at most
$$(d - a_i) + (f-e) - (d - e + 1) + (b_i - f) = b_i-a_i-1 = m -1.$$
This completes the proof of the lemma.
\end{proof}

\begin{thm}\label{thm:main}
Let $K$ be a field of arbitrary characteristic.  Then the ideal of
adjacent minors $I_{mn}(m)$ can be written as the irredundant
intersection of radical ideals
$$
I_{mn}(m) = \bigcap P_{\Gamma}
$$
where the intersection runs over all prime sequences of $[0,n+1]$.
When $char(K) = 0$  or when $m \leq 3$ in arbitrary
characteristic this is a minimal prime decomposition.
\end{thm}
\begin{proof}
Since $P_\Gamma$ is radical by Corollary
\ref{cor:rad}, the intersection
$\bigcap P_\Gamma$  is also radical.  Moreover, given any prime
sequence $\Gamma$, each adjacent $m \times m$ minor belongs to
$P_{\Gamma}$ since the column indices of this minor are either
contained in an interval $[a_i, b_i]$ in $\Gamma$ or they contain
the indices of one of the overlaps $[a_{i+1}, b_i]$.  This shows
that $I_{mn}(m)$ is contained in this radical ideal.  If $K$ is
algebraically closed, Lemma \ref{lem:algo} and the
Nullstellensatz imply that $I_{mn}(m)$ is equal to the
intersection.  Since all the ideals in question lie in
$K[x_{ij}]$ for any field $K$ we deduce that the equation 
holds over any field by passing to the algebraic closure.
In order to prove that this intersection is irredundant we need
to argue that if $\Gamma \neq \Gamma'$ then $P_\Gamma$ and
$P_{\Gamma'}$ are incomparable.  This is a consequence of our
Gr\"obner basis arguments and is proven in Corollary
\ref{cor:inc}. The intersection is a prime decomposition in
characteristic zero because $P_\Gamma$ is prime when $\mathrm{char}(K) 
= 0$:  this is the content of Theorem \ref{thm:prime}.  Similarly,
all the ideals $P_\Gamma$ are prime when $m \leq 3$ and the
characteristic is arbitrary.  This is proven in Corollary~\ref{cor:m3}.
\end{proof}

\begin{thm}\label{thm:counting}
Let  $f_m(n)$ be  the number of primes in the prime decomposition
of $I_{mn}(m)$. Then $f_m(n)$ is generated by the following
recurrence:
$$
f_m(n+1) = \sum_{i=0}^{m-1} f_m(n-i)
$$
subject to the initial conditions $f_m(1) = f_m(2) =
\cdots = f_m(m-2) = 0$, $f_m(m-1) = 1$ and $f_m(m)=1$.
\end{thm}

\begin{proof}
We count the prime sequences $\Gamma$ on $[0,n+1]$. There
are no such sequences when $n < m-1$ and there is a unique
sequence when $n = m-1$ or $n = m$. If the last interval
$[a_i, n+1]$ in $\Gamma$ has width greater than $m+1$ then
$\Gamma' = \Gamma - [a_i, n+1] \cup [a_i,n]$ is a prime
sequence of $[0,n]$. If the width of $[a_i, n+1]$ is $m+1$, then
$\Gamma' = \Gamma - [a_i, n+1]$ is a prime sequence of $[0,n+1-j]$
for $2 \leq j \leq m$. This gives an injective map from
the set of prime sequences of $[0,n+1]$ to the disjoint
union of prime sequences of $[0, n+1-m], [0,n+2-m], \ldots, [0,n]$.
It is also easy to see that the inverse of this map is injective.
Hence these two sets have the same cardinality which proves the
theorem.
\end{proof}

\section{A new class of prime determinantal ideals}

We now prove that the ideals $P_\Gamma$ are prime ideals in
characteristic zero.  We believe they are prime in arbitrary
characteristic and we  verify this conjecture in special cases.
First we will show that $P_\Gamma$ is a radical ideal
through a Gr\"obner basis argument which  does not
depend on $\mathrm{char}(K)$. Then we use an intricate
geometric argument to show that $\mathcal{V}(P_\Gamma)$ is
irreducible over fields of characteristic zero.

\subsection*{A Gr\"obner basis}
We will use the diagonal term order introduced in Proposition
\ref{complete-inter}. The argument will also depend on the
following lemma proved in \cite{Conca2}.

\begin{lemma} \label{lem:GB-lemma}
Let $I$ and $J$ be two homogeneous ideals of a polynomial ring
$K[x_1, \ldots, x_n]$, and let $F$ and $G$ be Gr\"obner bases
of $I$ and $J$ with respect to a fixed term order $\prec$. Then
$F \cup G$ is a Gr\"obner basis of $I+J$ with respect to
$\prec$ if and only if for every $f \in F$ and $g \in G$ there
exists $ h \in I \cap J$ such that $\mathrm{in}(h) =
\mathrm{LCM}(\mathrm{in}(f), \mathrm{in}(g))$.
\end{lemma}

\noindent Our main Gr\"obner basis result follows from the result below.

\begin{lemma} \label{lem:simple-GB}
Let $F$ be the set of $m \times m$ minors of $X_{mn}$ and let $G$
be the set of the $k \times k$ minors of the submatrix which
consists of either the first or the last
$k$ columns of $X_{mn}$ where $k < m$. Then with respect
to the lexicographic term
order $x_{11} \succ x_{12} \succ \cdots \succ x_{1n} \succ \cdots
\succ x_{mn}$ the set $F \cup G$ is a Gr\"obner basis of the ideal
it generates.
\end{lemma}
\begin{proof} We prove the case where $G$ is the set of
the $k \times k$ minors of the submatrix $Y$ consisting of the
first $k$ columns of $X_{mn}$ since the other case follows from a
symmetric argument similar to the one we give below. We will use
Lemma \ref{lem:GB-lemma} where $I = \langle F \rangle$ and $J =
\langle G \rangle$.  Note that $F$ and $G$ are Gr\"obner bases for
$I$ and $J$ with respect to the given term order by \cite{St2}.
For $f \in F$ and $g \in G$ we want to show that there is $h
\in I \cap J$ such that $\mathrm{in}(h) =
\mathrm{LCM}(\mathrm{in}(f),\mathrm{in}(g))$. We will construct
$h$ as follows:  let $\mathrm{in}(f) = x_{1i_1}x_{2i_2}\cdots
x_{mi_m}$ where $1 \leq i_1 < i_2 < \cdots < i_m \leq n$, and let
$\mathrm{in}(g) = x_{j_11}x_{j_22}\cdots x_{j_kk}$ where $1 \leq
j_1 < j_2 < \cdots < j_k \leq m$.  It is not hard to see that if
$\mathrm{in}(f)$ contains a variable $x_{si_s}$ where $i_s \leq k$
then for the corresponding variable $x_{j_{i_s}i_s}$ of
$\mathrm{in}(g)$ we have $j_{i_s} \geq s$.  Let $Y_1$ be the set
of columns of $Y$ indexed by the $j_t$ with $j_{i_s} = s$, and let
$Y_2$ be the set of those columns of $Y$ indexed by those $j_t$
which have $j_{i_s} > s$.  Moreover, let $Y_3$ be the set of
columns that do not contain a variable from $\mathrm{in}(f)$; that
is, $Y_3$ consists of the columns of $Y$ which are not in $Y_1$ or
$Y_2$.  Finally, $Y_4$ will be the set of columns of $X_{mn}$ with
indices $\{i_t: \, i_t
> k\}$. We make two simple observations. First of all, the sum
$|Y_1| + |Y_2| + |Y_4|$ is equal to $m$, and secondly, $Y_1$ comes
before all of the other $Y_i$ in $X_{mn}$: indeed, $Y_1$ is the
first $|Y_1|$ columns of $X_{mn}$. Now let us look at the rows of
$Y$ in which a variable of $\mathrm{in}(g)$ that is also either in
$Y_2$ or $Y_3$ appears. These rows form a $(|Y_2|+|Y_3|) \times k
$ submatrix of $Y$ that we  will denote by $A$.
With all this data we construct the $(m + |Y_2|+|Y_3|) \times (m +
|Y_2|+|Y_3|)$ matrix
$$\left[
\begin{array}{c|c}
A & 0 \\
\hline
Y & Y_2 \, | \,  Y_4
\end{array}
\right],
$$
and we let $h$ be its determinant. Since $h$ can be computed by
the Laplace expansion either using the $m \times m$ minors of the
last $m$ rows, or using the $k \times k$ minors of the first
$k$ columns we deduce that $h$ is in $I \cap J$.  The specific
term order we use together with the second observation above gives
us the fact that $\mathrm{in}(h) = \mathrm{LCM}(\mathrm{in}(f),
\mathrm{in}(g))$. This is the easiest to see by computing the Laplace 
expansion using the
first $|Y_2| + |Y_3|$ rows of the matrix.
\end{proof}

\begin{example} The proof of Lemma \ref{lem:simple-GB} relies on
the construction of a special element $h$ in $I \cap J$.  We will
now describe an example of this construction in the case $m = 5$,
$n = 6$, and $k = 3$ and we will suppose that we are taking $3
\times 3$ minors from the last three columns of $X_{mn}$.  In other
words, we illustrate the symmetrical case that we omitted in the 
above proof. We will
consider the special case where $f$ is the $5 \times 5$ minor with
column indices $\{1,2,3,4,6\}$ and $g$ is the $3 \times 3$ minor
with row indices $\{2,3,5\}$.  We can represent the situation
pictorially with a marked matrix:  the crosses $\times$ represent
variables which appear in the leading term of $f$ and the squares
$\Box$ represent variables which appear in the leading term of
$g$.  Our marked matrix is

$$\left[ \begin{array}{ccc|ccc}
\times &   &   &   &   &    \\
  & \times &   & \Box &  &  \\
  &   & \times &   & \Box &  \\
  &   &   & \times &  &   \\
  &   &   &   &  &  \! \boxtimes \end{array} \right]. $$

\noindent  According to the symmetric version of the construction,
we take $Y_1$ to consist of the last column of the matrix, $Y_2$ is
the third to last column, $Y_3$ is the second to last column, and
$Y_4$ consists of the first three columns.  We construct the new
matrix whose determinant is the desired polynomial $h$.  In this
new matrix, we again use symbols to mark the desired variables in
the leading term.  This new matrix is a $7 \times 7$ matrix and
looks like

$$\left[ \begin{array}{cccc|ccc}
0 & 0 & 0 & 0 & \Box &  &  \\
0 & 0 & 0 & 0 &   & \Box &  \\ \hline
\times &   &  &   &   &   &    \\
  & \times &    & \Box & \Box &  &  \\
  &   & \times & &   & \Box &  \\
  &   &   & \times  & \times &   &   \\
  &   &   &   &  &   &  \! \boxtimes \end{array} \right]. $$

\noindent  It is easy to see that $\IN(h) =
\mathrm{LCM}(\IN(f),\IN(g))$: just use the Laplace expansion
along the first two rows.
\end{example}

\begin{thm} \label{thm:GB}
With respect to the lexicographic term
order $x_{11} \succ x_{12} \succ \cdots \succ x_{1n} \succ \cdots
\succ x_{mn}$ all the minors defining $P_\Gamma$ form a
Gr\"obner basis.
\end{thm}
\begin{proof} We do induction on the number of intervals
in $\Gamma = \{[a_1,b_1], \ldots, [a_t,b_t]\}$. If
$\Gamma = \{[0,n+1]\}$, then $P_\Gamma$ is just generated
by the $m \times m$ minors of $X_{mn}$ and by
the results in \cite{St2} they form a Gr\"obner basis.
When there is more than one interval then $\Gamma' = \Gamma - [a_t, n+1]$
is a prime sequence for $[0,b_{t-1}+1]$. By induction, the set of
minors $F$ generating $P_{\Gamma'}$ is a Gr\"obner basis of $I :=
I_{mb_{t-1}}(m)$. Now we let $J$ be the ideal generated
by the $m \times m$ minors corresponding to the interval $[a_t,b_t]$
and the maximal minors of the overlap $[a_t, b_{t-1}]$.
We let $ k := b_{t-1} - a_t + 1$, and we denote the set of these
$k \times k$ minors  together with the $m \times m$ minors that generate
$J$ by $G$. Lemma \ref{lem:simple-GB} implies that $G$ is a
Gr\"obner basis of $J$. Now
we will use Lemma \ref{lem:GB-lemma} to prove the theorem.
Observe that if $f \in F$ and $g \in G$ are minors of
submatrices  corresponding to intervals or overlaps of intervals
which do not share a column, then
$\mathrm{LCM}(\mathrm{in}(f), \mathrm{in}(g)) = \mathrm{in}(f) \cdot
\mathrm{in}(g)$ and we choose $h = f \cdot g$. Hence we only need
to study the pairs of intervals that do overlap. Here is the list of the
cases we need to consider:
\begin{itemize}
\item[(a)] both $f$ and $g$ are $m \times m$ minors,
\item[(b)] $f$ is an $s \times s$ minor coming from an overlap
that also intersects the interval $[a_t,n+1]$, and $g$ is an
$m \times m$ minor,
\item[(c)] $f$ is as in (b), and $g$ is a $k \times k$ minor,
\item[(d)] $f$ is an $m \times m$ minor coming from an interval
that is not $[a_{t-1}, b_{t-1}]$ 
and $g$ is $k \times k$ minor, and
\item[(e)] $f$ is an $m \times m$ minor coming from $[a_{t-1}, b_{t-1}]$
and $g$ is a $k \times k$ minor.
\end{itemize}
The last case is covered by the proof of Lemma
\ref{lem:simple-GB}. In all the other cases, simple arguments show
that the leading terms of $f$ and $g$ are relatively prime and
hence we choose  $h = f \cdot g$.  For completeness, we go through
this argument for case (c).  The main tool is the following simple
observation.  For any maximal minor of any matrix, the leading
term selected by our diagonal lexicographic term order has all of
its variables lying in the parallelogram-shaped region bounded by
the diagonal extending from the upper left hand corner of the
matrix and the diagonal extending from the lower right hand
corner.  Since $\Gamma$ is a prime sequence, the smallest interval
$[a,b]$ which contains the column indices of both $f$ and 
$g$  has width greater than or
equal to $m +1$.  This ensures that the two regions corresponding
to the possible variables in the leading terms of these minors do
not intersect, because the diagonal from the upper left corner of
$X[a,b]$ is below the diagonal from the lower right corner of
$X[a,b]$.  This guarantees that the leading terms of $f$ and $g$
are relatively prime as desired.
\end{proof}

\begin{cor} \label{cor:rad} The ideal $P_\Gamma$ is radical.
\end{cor}
\begin{proof} The initial ideal of $P_\Gamma$ given
by Theorem \ref{thm:GB} is squarefree, and therefore it is radical.
Then $P_\Gamma$ is also radical.
\end{proof}

\begin{cor}\label{cor:inc}
If $\Gamma \neq \Gamma'$ then $P_\Gamma$ and $P_{\Gamma'}$ are
incomparable.
\end{cor}
\begin{proof}
We will show that $P_\Gamma$ is not contained in $P_{\Gamma'}$.
For this it suffices to show that there is a minor among the
generators of $P_\Gamma$ which is not contained in $P_{\Gamma'}$.
Let $[a_i,b_i]$ be the first interval of $\Gamma$ which is not
contained in $\Gamma'$ and let $[c_i,d_i]$ be the corresponding
$i$th interval of $\Gamma'$. The intervals $[a_1, b_1], \ldots
[a_{i-1},b_{i-1}]$ are the first $i-1$ intervals which are common
to both $\Gamma$ and $\Gamma'$. There are a few cases to consider.

If $i =1$ then $[a_1,b_1] = [0,b_1]$ and $[c_1, d_1] = [0,d_1]$.
Suppose that $b_1 > d_1$.  Among the indices in the interval
$[d_1+1, b_1]$ there exists at least one index $e$ 
so that $[e,e]$ is not 
an interval obtained by overlapping two consecutive intervals
in $\Gamma'$. Then
the $m \times m$ minor with columns indices $\{1, \ldots, m-1, e\}$ 
is contained in
$P_\Gamma$ but not in $P_{\Gamma'}$ because its leading term is
not divisible by any leading term in the Gr\"obner basis for
$P_{\Gamma'}$.  If we suppose that $b_1 < d_1$, then any $(b_1 -
a_2 +1) \times (b_1 - a_2 + 1)$ minor with column indices
$[a_2,b_1]$ belongs to $P_\Gamma$ but not $P_{\Gamma'}$ since its
leading term is not divisible by any leading term in the Gr\"obner
basis for $P_{\Gamma'}$.

Now we suppose that $i > 1$.  The arguments are similar to those
in the preceding paragraph and we sketch them briefly.  Suppose
$a_i < c_i$. Then there is an $m \times m$ minor with column
indices in $[a_i,b_i]$ using the column index $a_i$ which is
contained in $P_\Gamma$ but not $P_\Gamma'$.  If $a_i > c_i$ then
there is an $(b_{i-1} - a_i +1) \times (b_{i-1} - a_{i} + 1)$
minor with column indices equal to $[a_i, b_{i-1}]$ which is
contained in $P_{\Gamma}$ but not $P_{\Gamma'}$.  Finally, if $a_i
= c_i$ then a minor modification of the $i = 1$ case shows that
$P_\Gamma$ contains a minor which is not contained in
$P_{\Gamma'}$.
\end{proof}

\subsection*{$\mathcal{V}(P_\Gamma)$ is irreducible}

Before proceeding with the proof,  we will outline the strategy
that we will employ to show that $\mathcal{V}(P_\Gamma)$ is
irreducible over a field $K$ with $\mathrm{char}(K) = 0$. First,
we will construct a morphism from an irreducible variety
$\mathcal{X}$ to $\mathcal{V}(P_\Gamma)$. Then we will argue that
this morphism surjects onto a Zariski open subset $\mathcal{W}$ of
$\mathcal{V}(P_\Gamma)$ when restricted to a Zariski open (and
necessarily irreducible) subset $\mathcal{Y}$ of $\mathcal{X}$. This
implies that $\mathcal{W}$ is irreducible. Up to this point the
results will be obtained without any assumptions on the
characteristic of the field. Then we will assume that $K= \CC$,
and we will show that the closure of $\mathcal{W}$ is equal to
$\mathcal{V}(P_\Gamma)$ which proves that $\mathcal{V}(P_\Gamma)$
is irreducible. This will require a perturbation argument 
which we present in the next subsection. 
Finally, we use standard arguments in the proof of Theorem 
\ref{thm:prime} 
to show that $P_\Gamma$ is prime over any field of characteristic zero.

We first define the irreducible variety $\mathcal{X}$. In order to
do this we need to introduce a poset $\mathcal{Q}_\Gamma$
associated to a prime sequence $\Gamma$.

\begin{defn}
Let $\Gamma$ be a prime sequence. The elements of the poset
$\mathcal{Q}_\Gamma$ are certain subintervals of the intervals
in $\Gamma$ which will be defined recursively, and these subintervals
are ordered with respect to inclusion.
The intervals in $\Gamma$ are the maximal elements
of $\mathcal{Q}_\Gamma$, and we sort them with respect
to each interval's starting index, the {\em left border},
in ascending order. These will form the elements
in {\em row} $1$. The elements in row $2$ are 
the nonempty subintervals obtained
by intersecting two consecutive intervals in row $1$.
We also sort row $2$ in ascending order with respect to the left borders.
The subsequent rows are defined recursively: the elements in
row $r$ consist of all nonempty intervals that arise from the intersection
of two consecutive elements from row $r-1$.
It is clear that every nonmaximal element is covered
by exactly two elements (a left and a right parent), and
each nonminimal element covers at most two other elements
(a left and a right child).
\end{defn}

\begin{figure} \label{fig:poset}
\centerline{ \epsfig{file=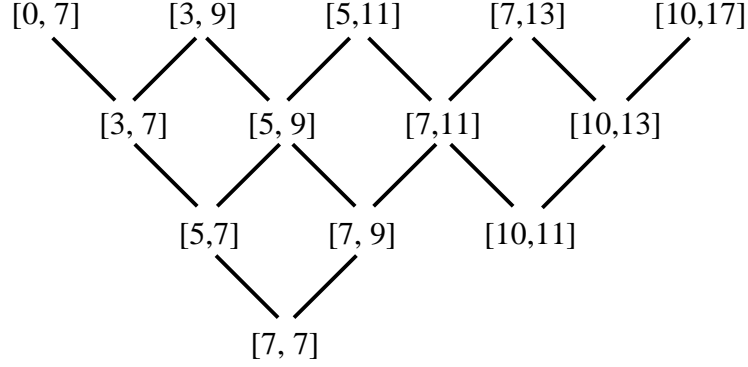}} \caption{The poset
$\mathcal{Q}_\Gamma$ for $\Gamma = \{[0,7],[3,9],[5,11],[7,13],
[10,17] \}$}
\end{figure}

\begin{example}
Let $m = 6$ and consider the sequence of intervals
$$ \Gamma \,\, = \,\, \{[0,7],[3,9],[5,11],[7,13], [10,17] \}.$$
The second row of the poset
consists of the overlapping intervals $[3,7],[5,9], [7,11]$,
and $[10,13] $. The third row is formed by  the intervals
$[5,7], [7,9]$, and $[10,11]$.  The fourth and final row of the poset
is the interval $[7,7]$. This poset is
illustrated in Figure \ref{fig:poset}.
\end{example}

In order to define $\mathcal{X}$ we need one more piece of
information. This will be a positive integer attached to each
element of $\mathcal{Q}_\Gamma$.

\begin{defn}
For each $p \in \mathcal{Q}_\Gamma$ let

$$ D(p) := \left\{  \begin{array}{cl}
m-1 & \mbox{if $p$ is in the first row of  } \mathcal{Q}_\Gamma \\
w(p) -1 & \mbox{ if $p$ is in the second row of  }  \mathcal{Q}_\Gamma  \\
w(p)  & \mbox{otherwise} \end{array} \right.
$$
\noindent where $w(p)$ is the width of the interval $p$.
\end{defn}

Now each element $p \in \mathcal{Q}_\Gamma$ will give rise to a
general linear group $\gl_{k(p)}$ of invertible $k(p) \times k(p)$
matrices where $k(p) = D(q)-D(q')$, and $q$ is the left
parent of $p$ and $q'$ is the left child of $p$. If $p$ does not
have a left parent then we set $D(q) = m$, and if $p$ does not
have a left child, then we set $D(q') = 0$. Moreover each maximal
element $q \in \mathcal{Q}_\Gamma$ will give rise to an affine
space $\AA^{\ell(q)}$, and we define $\ell(q)$ as follows: suppose
$q$ corresponds to the interval $[a_s, b_s]$ and let $[a_{s+1},
b_{s+1}]$ be the next interval (if there is one). Let $\Lambda :=
[a_s, a_{s+1}-1]$ or $\Lambda := [a_s, n]$ if $[a_s, b_s]$ is the
last interval. Now for each index $i \in \Lambda$ there is a
unique $p(i) \in \mathcal{Q}_\Gamma$ which is minimal among all
elements containing $i$. It is an easy exercise to see that $p(i)
\in \{p_0, \ldots, p_r\}$ where $p_0 =q $ and $p_{j+1}$
is the left child of $p_j$. With this we define $\ell(q) = \sum_{i
\in \Lambda} D(p(i))$. Finally we arrive at the variety
$$ \mathcal{X} \,\, := \,\, \prod_{p \in \mathcal{Q}_\Gamma} \gl_{k(p)}
\,\,\,  \times \,\,\, \prod_{q \,\, maximal} \AA^{\ell(q)}. $$

\noindent We note that over an infinite field $\mathcal{X}$
is irreducible since it is the product of irreducible varieties.

Next we define a map $\phi$ from $\mathcal{X}$ to $\AA^{mn}$, the
space of all $m \times n$ matrices. Given a point in
$x \in \mathcal{X}$
we will build an $m \times n$ matrix piece by piece using
the intervals in $\Gamma$. We start with the last interval $[a_t, n+1]$
and the corresponding maximal element $q \in \mathcal{Q}_\Gamma$.
Then we build an $m \times |\Lambda|$ matrix $\mathbf{Z}$ as follows: 
for each $i \in \Lambda$ we set all entries
in column $i$ with row indices $D(p(i))+1, D(p(i))+2, \ldots, m$
to zero. There are precisely $\ell(q)$ entries in $\mathbf{Z}$ that are
not set to zero yet, and we ``plug in'' the coordinates of
the point $x$ corresponding to $\AA^{\ell(q)}$ to these entries.
We set $\mathbf{X} := \mathbf{Z}$. 
Now let $q=p_0, p_1, \ldots, p_s$ be the elements of $\mathcal{Q}_\Gamma$
such that $p_{j+1}$ is the left child of $p_j$,
and let $g_j \in \gl_{k(p_j)}$ be the matrices that could be read off
from the corresponding coordinates of $x$. For $j =0, \ldots, s$
we define $\mathbf{X} := g_j \cdot \mathbf{X}$ recursively, where
$g_j \cdot \mathbf{X}$ is obtained by multiplying the last
$k(p_j)$ rows of the first $D(q_j)$ rows of
$\mathbf{X}$, and $q_j$ is the left parent of $p_j$
(since $D(q_j) \geq k(p_j)$, by the definition of $k(p_j)$ this makes
sense).

After we have gone through the sequence $p_0, \ldots, p_s$, let the
resulting matrix be $\mathbf{Y}$.
Next we move onto the second to last interval $[a_{t-1}, b_{t-1}]$, and
using the set $\Lambda$ associated to this interval we build
a matrix $\mathbf{Z}$, and then we set
$\mathbf{X} := [\mathbf{Z} | \mathbf{Y}]$. Now using the various
invertible matrices associated to the sequence of the left children 
starting from 
$[a_{t-1}, b_{t-1}]$ we repeat this procedure. 
Clearly the result of this construction is an $m \times n$ matrix.
It is also clear that this map is a polynomial map and hence
a morphism.

\begin{example} \label{ex:detail}
This is a detailed example displaying
the variety $\mathcal{X}$ and the recursive construction of the map
$\phi$.  Let $m = 4$ and let $\Gamma $ be the prime sequence
$\Gamma = \{[0,5], [3,7], [5,10] \}$.  The second row of the poset
$\mathcal{Q}_\Gamma$ consists of the two intervals $[3,5]$ and
$[5,7]$, and the third row of the poset is the singleton interval
$[5,5]$.  According to the construction of $\mathcal{X}$ we have

$$ \mathcal{X} = \gl_4 \times \gl_3 \times \gl_2 \times \gl_2
\times \gl_2 \times \gl_2 \times \AA^6 \times \AA^4 \times
\AA^{11}.
$$

\noindent  We have ordered the general linear groups and the
affine spaces in the \emph{reverse} of the order in which they are
used in the map $\phi$.  This should not be confusing to the
reader:  the ordering of the general linear groups mimics the
right to left order of group actions and the affine spaces are
ordered in this way as a reminder that we construct the matrix in the
image of $\phi$ from right to left.
Now let $x$ be an arbitrary point in the variety $\mathcal{X}$.
We begin with the interval
$[5,10]$, the last interval in $\Gamma$, and use the affine space
$\AA^{11}$ to construct a $4 \times 5$ matrix $\mathbf{Z}$ which
looks like
$$ \mathbf{Z} = \left[  \begin{array}{ccccc}
* & * & * & * & * \\
0 & * & * & * & * \\
0 & 0 & 0 & * & * \\
0 & 0 & 0 & 0 & 0 \end{array} \right]$$

\noindent and corresponds to columns 5 through 9 of our eventual
completed matrix.  We set $\mathbf{X} := \mathbf{Z}$. 
Now we read down the right-most chain in the
poset and apply the action of general linear groups accordingly.
In particular, we apply $g_1 \in \gl_2$ to the bottom two rows of
$\mathbf{X}$, then $g_2 \in \gl_2$ to the middle two rows of $g_1
\cdot \mathbf{X}$, and finally $g_3 \in \gl_2$ to the first two rows 
of $g_2 \cdot g_1
\cdot \mathbf{X}$.  Pictorially, we have
\fontsize{10}{12}
$$
\left[  \begin{array}{ccccc}
* & * & * & * & * \\
0 & * & * & * & * \\
0 & 0 & 0 & * & * \\
0 & 0 & 0 & 0 & 0 \end{array} \right] {g_1 \atop \longrightarrow}
\left[  \begin{array}{ccccc}
* & * & * & * & * \\
0 & * & * & * & * \\
0 & 0 & 0 & * & * \\
0 & 0 & 0 & * & * \end{array} \right] {g_2 \atop \longrightarrow}
\left[  \begin{array}{ccccc}
* & * & * & * & * \\
0 & * & * & * & * \\
0 & * & * & * & * \\
0 & 0 & 0 & * & * \end{array} \right] {g_3 \atop \longrightarrow}
\left[  \begin{array}{ccccc}
* & * & * & * & * \\
** & * & * & * & * \\
0 & * & * & * & * \\
0 & 0 & 0 & * & * \end{array} \right]. $$

\normalsize

\noindent  where the last matrix is the matrix $\mathbf{Y}$ obtained
at the end of this iteration of the construction. Now we look at
the second to last interval $[3,7]$ in $\Gamma$.
Comparing with the interval $[5,10]$ we
see that $\Lambda = [3,4]$, and we add two new columns $\mathbf{Z}$ 
to $\mathbf{Y}$ above. 
These come from our $\AA^4$ to arrive at a matrix $\mathbf{X} :=
\left[\mathbf{Z} | \mathbf{Y}\right]$. Reading the second descending
chain in $\mathcal{Q}_{\Gamma}$ we apply $g_4 \in \gl_2$ to
the last two rows of $\mathbf{X}$ and then apply $g_5 \in \gl_3$
to the first three rows of $g_4 \cdot \mathbf{X}$.  Pictorially,
this looks like
\fontsize{10}{12}
$$
\left[ \begin{array}{cc|ccccc}
*  & * & * & * & * & * & * \\
** & * & * & * & * & * & * \\
0  & 0 & 0 & * & * & * & * \\
0  & 0 & 0 & 0 & 0 & * & * \end{array} \right] {g_4 \atop
\longrightarrow} \left[ \begin{array}{cc|ccccc}
*  & * & * & * & * & * & * \\
** & * & * & * & * & * & * \\
0  & 0 & 0 & * & * & * & * \\
0  & 0 & 0 & * & * & * & * \end{array} \right] {g_5 \atop
\longrightarrow}
\left[ \begin{array}{cc|ccccc}
*  & * & * & * & * & * & * \\
** & * & * & * & * & * & * \\
** & * & * & * & * & * & * \\
0  & 0 & 0 & * & * & * & * \end{array} \right].
$$

\normalsize
\noindent  And again the last matrix is the matrix $\mathbf{Y}$ obtained
at the end of the second iteration of the construction.
We are now at the last step and we adjoin two new columns $\mathbf{Z}$
to our matrix $\mathbf{Y}$. The entries in these columns
come from the $\AA^6$. We form  the matrix
$$\mathbf{X} := \left[\mathbf{Z} | \mathbf{Y}\right]
= \left[ \begin{array}{cc|ccccccc}
* & * & *  & * & * & * & * & * & * \\
** & * & *  & * & * & * & * & * & * \\
** & * & *  & * & * & * & * & * & * \\
0 & 0 & 0  & 0 & 0 & * & * & * & * \end{array} \right]. $$
 Since the
interval $[0,5]$ has no left child, we deduce that
we should apply $g_6 \in \gl_4$ to the entire matrix.  This final
matrix $g_6 \cdot \mathbf{X}$ is the image of $x$ under $\phi$.
\end{example}

\begin{prop} The image of $\phi$ is contained in $\mathcal{V}(P_\Gamma)$.
\end{prop}
\begin{proof} We need to show that for every $x \in \mathcal{X}$
all the minors that generate $P_\Gamma$ vanish on $\phi(x)$.
We will prove this by using the definition of $\phi$. First we
observe that if a set of minors vanish on the partial matrix
$\mathbf{X}$ in the definition of $\phi$ then after the row operations
$g_j \cdot \mathbf{X}$ these minors will still vanish on $\mathbf{X}$.
We will show that as we build $\mathbf{X}$ each submatrix
of $\mathbf{X}$ that corresponds to $p \in \mathcal{Q}_\Gamma$
has rank at most $D(p)$. This is certainly true after constructing
$\mathbf{X}$ corresponding to the last interval of $\Gamma$, since
at most the first $D(p_j)$ rows of each submatrix corresponding to
$p_j$ are nonzero. An inductive argument shows that after applying
$g_j$ to $\mathbf{X}$, the columns of $\mathbf{X}$ that are
in the submatrix corresponding to $p_k$ for $k=0, \ldots, j$, but that
are not in the submatrix corresponding to $p_{k+1}$ have
nonzero elements in at most the first $D(q_k)$ rows where $q_k$ is the
left parent of $p_k$. So when
$\mathbf{Y}$ is constructed at most the first $D(q_j)$
rows of the matrix corresponding to $p_j$ for $j=0, \ldots, s$ are
nonzero.
In order to finish the proof
by induction, we assume that after constructing the matrix $\mathbf{Y}$
for an interval in $[a_r, b_r]$ where $r >1$ all the minors
arising from the intervals $[a_r, b_r], [a_{r+1}, b_{r+1}], \ldots,
[a_t, n+1]$ and their consecutive overlaps vanish on $\mathbf{Y}$, and
in the submatrices corresponding to the sequence $p_0, \ldots, p_s$
(where $p_0$ is $[a_r,b_r]$) at most the first $D(q_j)$
rows are nonzero, where $q_j$ is the left parent of $p_j$.
When we move to the next
interval $[a_{r-1}, b_{r-1}]$ with the corresponding sequence
of elements $\bar{p}_0, \ldots, \bar{p}_u$, first we consruct
$\left[\mathbf{Z} \, | \mathbf{Y}\right]$. It is easy to see that
the submatrix $\mathbf{A}_j$ of this matrix corresponding to $\bar{p}_j$
is obtained
by concatenating the portion of $\mathbf{Z}$
contained in  $\mathbf{A}_j$ with the submatrix corresponding to
the right child of $\bar{p}_j$. Now at most the first $D(\bar{p}_j)$
rows of the portion of $\mathbf{A}_j$ contained in $\mathbf{Z}$
are nonzero, and by induction the same is true for the submatrix
corresponding to the right child of  $\bar{p}_j$. Hence at most
the first $D(\bar{p}_j)$ rows of $\mathbf{A}_j$ are nonzero. This
shows that the minors arising from the intervals $[a_{r-1}, b_{r-1}],
[a_r, b_r], [a_{r+1}, b_{r+1}], \ldots, [a_t, n+1]$
and their consecutive overlaps vanish on $\mathbf{X} :=
\left[\mathbf{Z} \, | \mathbf{Y}\right]$. 
After applying the row operations $\bar{g}_j$, at most
the first $D(\bar{q}_j)$ rows
of the matrix corresponding to $\bar{p}_j$ will be nonzero
where $\bar{q}_j$ is the left parent of $\bar{p}_j$ because 
$k(\bar{p}_j) = D(\bar{q}_j) - D(\bar{p}_{j+1})$. This implies that 
$\mathbf{X}$ has the properties the induction is based on, and 
this completes the induction.
\end{proof}

Now we let $\mathcal{W}$ be the subset of $\mathcal{V}(P_\Gamma)$
consisting of matrices $\mathbf{X}$ where the rank of  each
submatrix of $\mathbf{X}$ corresponding to $p \in
\mathcal{Q}_\Gamma$ is equal to $D(p)$. Since this subset is
defined by the non-vanishing of certain minors we conclude that it
is a Zariski open subset of $\mathcal{V}(P_\Gamma)$.  It is guaranteed to
be nonempty by the results in the next subsection.
Moreover, we
let $\mathcal{Y}$ be the set of $x \in \mathcal{X}$ such that
$\phi(x) \in \mathcal{W}$. We argue that $\mathcal{Y}$ is
an open subset of $\mathcal{X}$.  For this, consider an $x \in
\mathcal{X}$ where we take the entries as indeterminates. Then
$\phi(x)$ is a matrix with polynomial entries in the coordinates
of $x$.  Thus $\mathcal{Y}$ is defined by the non-vanishing of
certain minors of $\phi(x)$.
Furthermore, $\mathcal{Y}$ is irreducible since $\mathcal{X}$ is
irreducible.

\begin{prop} \label{prop:irred} The morphism $\phi \, : \, \mathcal{Y} \longrightarrow
\mathcal{W}$ is surjective, and therefore $\mathcal{W}$ is
irreducible when $K$ is an infinite field.
\end{prop}
\begin{proof} Since the second statement follows from the
first we just prove the first claim. We will do this
by constructing $x \in \mathcal{Y}$ for each
$\mathbf{X} \in \mathcal{W}$ such that $\phi(x) = \mathbf{X}$.
We start with the first interval
$p = [a_1, b_1]$ in $\Gamma$.
Since $\mathbf{X}[a_1,b_1]$ has rank $D(p) = m-1$, we can
find a $g \in \gl_{k(p)}$ where $k(p) = m$ so that $g \cdot
\mathbf{X}[a_1, b_1]$ is row-reduced, in particular, the last row
is a zero row.
We let $\mathbf{X} = g \cdot \mathbf{X}$, and we record
$g^{-1}$
as well as the entries of the first $D(p) = m -1$ rows of each column
of $\mathbf{X}$ with column index $i \in \Lambda$ (see the definition
of $\Lambda$ in the paragraph before Example \ref{ex:detail})
as part of the element $x$ we are constructing.
Then we delete these columns from $\mathbf{X}$ to obtain
the new $\mathbf{X}$.

By induction, suppose we have gone through the intervals
$[a_1,b_1], \ldots, [a_{r-1}, b_{r-1}]$ and the matrix $\mathbf{X}
= \mathbf{X}[a_r, n]$ has the following properties. Let $p_0, p_1,
\ldots, p_s$ be the sequence of elements where $p_0$ is $[a_{r-1},
b_{r-1}]$ and $p_{j+1}$ is the left child of $p_j$. Then the only
nonzero rows of the submatrix of $\mathbf{X}$ corresponding to the 
right child of $p_j$ are the first $D(p_j)$ rows of this
submatrix.

Now we let $\bar{p}_0, \ldots, \bar{p}_u$ be the sequence where
$\bar{p}_0$ is the interval $[a_{r}, b_{r}]$ and $\bar{p}_{j+1}$
is the left child of $\bar{p}_j$. We observe that, by induction,
only the first $D(\bar{q}_j)$ rows of the submatrix of
$\mathbf{X}$ corresponding to $\bar{p}_j$ could be nonzero where
$\bar{q}_j$ is the left parent of  $\bar{p}_j$. Now we apply
$\bar{g}_j \in \gl_{k(\bar{p}_j)}$ to $\mathbf{X}$ successively
starting from $j = u$ and finishing with $j=0$. In this process
$\bar{g}_j$ will be chosen as the matrix which will be applied to
the last $k(\bar{p}_j)$ rows of the first $D(\bar{q}_j)$ rows of
$\mathbf{X}$ so that after applying $\bar{g}_j$, the submatrix of
$\mathbf{X}$ corresponding to $\bar{p}_j$ is row reduced.  Since
we have assumed that the submatrix corresponding to
$\bar{p}_{j+1}$ has rank equal to $D(\bar{p}_{j+1})$, this implies
that the submatrix corresponding to $\bar{p}_j$ with row indices
$D(\bar{p}_{j+1}) + 1 $, $D(\bar{p}_{j+1}) + 2$, $\ldots$,
$D(\bar{p}_{j+1}) + k(\bar{p}_j)$ has rank $k(\bar{p}_j)-1$. 
This implies that after applying $\bar{g}_j$
the rows indexed by $D(\bar{p}_j)+1, D(\bar{p}_j)+2, \ldots, m$ in
the submatrix corresponding to $ \bar{p}_j $ will consist of
zeros.  Moreover, when applying $\bar{g}_j$, the definition of
$k(\bar{p}_j)$ and the particular rows which will be affected
guarantee that the zero rows of the submatrices corresponding to
$\bar{p}_{j+1}, \ldots, \bar{p}_u$ stay as zero rows. Hence when
we compute $\mathbf{X} := \bar{g}_0 \cdot \mathbf{X}$, we return
to the property we started with at the beginning of the induction
step, namely: the only nonzero rows of the submatrix of
$\mathbf{X}$ corresponding to $\bar{p}_j$ are the first
$D(\bar{p}_j)$ rows of this submatrix.

Now we delete the submatrix $\mathbf{X}[a_r, a_{r+1}-1]$ from
$\mathbf{X}$ to obtain the new $\mathbf{X}$ for the next
iteration, and we record the $\ell(\bar{p}_0)$ possibly nonzero
elements in the deleted columns of $\mathbf{X}$ as part of $x$
(this belongs to $\AA^{\ell(\bar{p}_0)}$) as well as the inverses
of all the matrices $\bar{g}_{j} \in \gl_{k(p_j)}$ which we used.
Since we have returned $\mathbf{X}$ to the form of the inductive
hypothesis, this shows that we can continue the procedure to
compute an $x \in \mathcal{Y}$ whose image under $\phi$ is
$\mathbf{X}$.

\end{proof}

We conclude this section with the proof that in certain special
cases, the map $\phi: \mathcal{X} \rightarrow
\mathcal{V}(P_\Gamma)$ is, in fact, surjective (i.e. not just
surjective on an open subset).  Hence, in these cases we may
conclude that $P_\Gamma$ is prime without resorting to the
analytic techniques in Proposition \ref{prop:perturb}.

\begin{prop}\label{prop:2overlap}
Suppose that $\mathcal{Q}_{\Gamma}$ has only two rows.  Then the
map $\phi: \mathcal{X} \rightarrow \mathcal{V}(P_\Gamma)$ is
surjective.
\end{prop}

\begin{proof}
We will closely follow the proof of Proposition \ref{prop:irred}
but with an extra twist.  Given an $\mathbf{X} \in
\mathcal{V}(P_\Gamma)$ we will construct $x \in \mathcal{X}$ such
that $\phi(x) = \mathbf{X}$.  We start with the first interval
$p = [a_1, b_1]$ in $\Gamma$. Since $\mathbf{X}[a_1,b_1]$ has
rank $D(p) = m-1$, we can find a $g \in \gl_{k(p)}$ where $k(p) =
m$ so that $g \cdot \mathbf{X}[a_1, b_1]$ is row-reduced, in
particular, the last row is a zero row. We let $\mathbf{X} = g
\cdot \mathbf{X}$, and we record $g^{-1}$ as well as the entries
of the first $D(p) = m -1$ rows of each column of $\mathbf{X}$
with column index $i \in \Lambda$ as part of the element $x$ we
are constructing. Then we delete these columns from $\mathbf{X}$
to obtain the new $\mathbf{X}$.

Since $\mathcal{Q}_\Gamma$ has only two rows, our inductive
hypothesis is simpler than Proposition \ref{prop:irred}.  Namely,
suppose that we have gone through the intervals $[a_1,b_1],
\ldots, [a_{r-1}, b_{r-1}]$ and the submatrix of $\mathbf{X} =
\mathbf{X}[a_r, n]$ indexed by $\Lambda = [a_r, b_{r-1}]$ has its
last row as a zero row.

We will let $\bar{p}_0 = [a_r, b_r]$ and  $\bar{p}_1 = [a_r, b_{r-1}]$.  
Since
$\mathcal{Q}_\Gamma$ has only two rows, these are the only two
elements of $\mathcal{Q}_\Gamma$ which we need to consider when we
perform our induction.  First we use an element of $g_1 \in
\gl_{m-1}$ to row reduce the submatrix consisting of the first
$k(\bar{p}_1) = m-1$ 
rows of $\mathbf{X}$ and the columns indexed by $\bar{p}_1$. This
submatrix has rank $\leq D(\bar{p}_1) = b_{r-1} - a_r$: if the rank of the
submatrix is strictly less than $D(\bar{p}_1)$ we must perform our row
reductions with caution to ensure that the submatrix of
$\mathbf{X}$ with columns indexed by $[b_{r-1} +1, b_r]$  and
consisting of the last $k(\bar{p}_0) = m - D(\bar{p}_1)$ 
rows has rank less than $m -
D(\bar{p}_1)$.  To ensure this possibility, we note that there are two
cases to consider.  In the first case,  the submatrix consisting of 
its first $m-1$ rows of $\mathbf{X}[a_r,b_r]$ has rank $m-1$.
In this case we can choose $g_1$ so that the $(m-1)$st row of $g_1
\cdot \mathbf{X}[a_r,b_r]$ is a multiple of the last row of
$\mathbf{X}[a_r,b_r]$.  In the second case, 
the submatrix consisting of the first $m-1$ rows of $\mathbf{X}[a_r,b_r]$
rank $< m-1$.  Then we can choose $g_1$ so that the $(m-1)$st row of $g_1
\cdot \mathbf{X}[a_r,b_r]$ is the zero row.  In either case,  this
ensures that the last $m-D(\bar{p}_1)$ rows of $\mathbf{X}[b_{r-1} +1,
b_r]$ has rank less than $m - D(\bar{p}_1)$.  Now we apply row reduction
via  $g_0$ in $\gl_{k(\bar{p}_0)}$ to the last $m - D(\bar{p}_1)$ rows of
$\mathbf{X}$ to bring $\mathbf{X}$ into the form
of the inductive hypothesis.

To complete the proof, we record the entries in first $D(\bar{p}_1)$ rows
of $\mathbf{X}[a_{r}, b_{r-1}]$, and 
the entries in the first  $D(\bar{p}_0)$ rows of 
$\mathbf{X}[b_{r-1}+1, a_{r+1}-1]$ (this becomes a set of
entries in $\AA^{\ell(\bar{p}_0)}$). 
We also record the inverses of $g_1$ and
$g_0$, and we delete the first $|\Lambda| = a_{r+1}-a_r$ 
columns from $\mathbf{X}$
to arrive at $\mathbf{X} := \mathbf{X}[a_{r+1},n]$.  By our
construction, this matrix is in proper form of the inductive
hypothesis, and so we may continue the process to construct $x$
such that $\phi(x) = \mathbf{X}$.
\end{proof}

\begin{cor} \label{cor:2overlap}
Let $K$ be a field of arbitrary characteristic and suppose that
$\mathcal{Q}_{\Gamma}$ has only two rows.  Then $P_\Gamma$ is a
prime ideal.
\end{cor}
\begin{proof}
If $K$ is algebraically closed, Proposition \ref{prop:2overlap} and
Corollary \ref{cor:rad} together with the Nullstellensatz imply
that $P_\Gamma$ is a prime ideal.  But this implies $P_\Gamma$ is
prime over any field by passing to the algebraic closure.
\end{proof}

\begin{cor}\label{cor:m3}
If $m \leq 3$, then $P_\Gamma$ is prime for any prime sequence
$\Gamma$.
\end{cor}
\begin{proof}
For $m \leq 2$ the statement was proven in \cite{DES}. When $m =
3$, each interval $[a_i,b_i] \in \Gamma$ has width greater than
or equal to 4, whereas each of the overlapping intervals
$[a_i,b_{i-1}]$ and $[a_{i+1},b_i]$ has width less than or equal
to 2.  This implies that the intervals $[a_i,b_{i-1}]$ and
$[a_{i+1},b_i]$ do not overlap and so $\mathcal{Q}_\Gamma$ has
only two rows.  By Corollary \ref{cor:2overlap}, $P_\Gamma$ is a
prime ideal.
\end{proof}

The reader may wonder why we have not shown that $\phi:\mathcal{X}
\rightarrow \mathcal{V}(P_\Gamma)$ is surjective in general,
eliminating the need for the analytic arguments in Proposition
\ref{prop:perturb}.  In general, it is not clear if this is true;
so we state it as a question.

\begin{question}
Is the morphism $\phi:\mathcal{X} \rightarrow
\mathcal{V}(P_\Gamma)$ always surjective?
\end{question}

We do not even know the answer in the case $m = 4$ with $\Gamma =
\{[0,5], [3,7], [5,10] \}$ from Example \ref{ex:detail}, which is
essentially the smallest instance not covered by Proposition
\ref{prop:2overlap}.  An affirmative answer to this question would
imply that $P_\Gamma$ is prime for all $\Gamma$ and in arbitrary
characteristic.

\subsection*{The perturbation argument}

We now present the details of the argument that every point of
$\mathcal{V}(P_{\Gamma})$ is arbitrarily close to $\mathcal{W}$
when the underlying field is $\CC$.  It suffices to show that
given a matrix $\mathbf{X} \in 
\mathcal{V}(P_{\Gamma}) \setminus \mathcal{W}$, 
there exists an infinitesimal perturbation which will
make the rank of all the submatrices corresponding to $q \in
\mathcal{Q}_\Gamma$ equal to $D(q)$.  Making this
perturbation requires care, since an arbitrary perturbation might
force the rank of some submatrix to jump to a value greater than
$D(q)$, and this will result in a matrix that is no longer in
$\mathcal{V}(P_{\Gamma})$.

For notational convenience we denote by
$\mathcal{Q}_\Gamma(\mathbf{X})$ the poset $\mathcal{Q}_\Gamma$
where the elements are taken to be the actual submatrices instead
of the intervals. This way, for instance, we will be able to work
with $\mathrm{span}(p)$ of $p \in \mathcal{Q}_\Gamma(\mathbf{X})$
which will mean the vector space spanned by the columns of $p$.
Similarly $\dim(p)$ will denote the dimension of this vector space.

\begin{defn}
Let $p$ be an element of the poset $\mathcal{Q}_\Gamma(\mathbf{X})$.
We let $\mathcal{M}(p)$ be the set of elements of
$\mathcal{Q}_\Gamma(\mathbf{X})$ above $p$
whose rank is equal to the desired maximal rank:

$$ \mathcal{M}(p) = \left\{ q \in  \mathcal{Q}_\Gamma(\mathbf{X}) |
q \geq p, \,\, \dim(q) = D(q)
\right\}.$$

\noindent Next we define
$Per(p)$, the \emph{vector space of allowable perturbations} to be

$$Per(p) = \left\{ \begin{array}{cl}
\bigcap_{q \in \mathcal{M}(p)} \mathrm{span}(q) &  \mbox{if  }
\mathcal{M}(p) \neq
\emptyset \\
\mathbb{C}^m & \mbox{if  } \mathcal{M}(p) = \emptyset \end{array}
\right. $$
\end{defn}

\begin{lemma}\label{lem:enough}
Let $p$ be an element of the poset $\mathcal{Q}_\Gamma(\mathbf{X})$. Then

$$\label{eq:enough} D(p) \leq \dim Per(p),$$

\noindent that is, there is a large enough vector space in which
perturbations can be made.
\end{lemma}

\begin{proof}
We suppose throughout that $\dim(p) < D(p)$ since in
the case of equality there is nothing to prove.  Assuming this, the case
where $\mathcal{M}(p) = \emptyset$ is trivial, so suppose
$\mathcal{M}(p)$  is nonempty.  If $p$ is in the first row
of $\mathcal{Q}_\Gamma(\mathbf{X})$ there is
nothing to show. If $p$ is in the second row then $\dim Per(p)
\geq m-2$ whereas $D(p) \leq m-2$ by the definition of $\Gamma$.
So suppose that $p$ is in
at least the third  row of the poset.

Clearly it is enough to take the minimal elements
of $\mathcal{M}(p)$ when computing $Per(p)$.  Furthermore, since
$\dim(p) < D(p) = w(p)$ we see that no $q > p$ can
belong to $\mathcal{M}(p)$ if $q$ is in at least the third row. 
Otherwise for
such a $q$ to be in $\mathcal{M}(p)$ would require that $q$ has
its rank equal to its width.  But then $\dim(p) = w(p)$
which is a contradiction. With this in mind, we first prove
the inequality in the statement of the lemma when the minimal elements
of $\mathcal{M}(p)$ consist of elements in the
second row of the poset.

Let $q_1, q_2, \ldots, q_r$ be all the elements in the second row
of the poset which are larger than $p$.  We assume that each $q_j$
spans a subspace of dimension $w(q_j) -1 := m - i_j - 1$
where $i_j > 0$.  Now
consider the intersection of the vector spaces spanned by the
$q_i$.  Since $q_1$ and $q_2$ are submatrices of $q_1 \vee q_2$
which has rank less than or equal to $m-1$, we deduce that the
vector space $\mathrm{span}(q_1) \cap \mathrm{span}(q_2)$ has
dimension at least $m - i_1 - i_2 - 1$.  By induction, the vector
space

$$ Per(p) = \bigcap_{j=1}^r \mathrm{span}(q_j)$$

\noindent has dimension $\geq m -1 - \sum_j i_j.$  On the other
hand, $w(p) \leq m - r + 1 - \sum_j i_j$ since the width of the
intervals in $\Gamma$ is at least $m+1$, and this completes the
proof in the case when the minimal elements of $\mathcal{M}(p)$
consist of elements in the second row of the poset.  The general
case now follows because removing one of the $q_j$ from
$\mathcal{M}(p)$ (and possibly adding something from the first
row)  can only make $\dim Per(p)$ larger.
\end{proof}

Now we show how
perturbations should be made inside a given rank-deficient matrix
$\mathbf{X}$
so that every submatrix corresponding to
$p \in \mathcal{Q}_\Gamma(\mathbf{X})$ has maximal rank $D(p)$.

\begin{prop} \label{prop:perturb}
Let $\mathbf{X} \in \mathcal{V}(P_\Gamma)$ be a matrix such that
$\dim(p) < D(p)$ for an element $p \in
\mathcal{Q}_\Gamma(\mathbf{X})$. Then there is an infinitesimal
perturbation of $\mathbf{X}$ to $\mathbf{X}' \in
\mathcal{V}(P_\Gamma)$ such that the rank of the corresponding
$p'$ in $\mathcal{Q}_\Gamma(\mathbf{X})$ increases, and $\dim(q)$
for any other element in $\mathcal{Q}_\Gamma(\mathbf{X})$ does not
decrease.
\end{prop}

\begin{proof}
We can assume that $p$ is minimal in
$\mathcal{Q}_\Gamma(\mathbf{X})$
among the submatrices that are rank-deficient.
It suffices to show that we can increase the rank of this
submatrix by one.  There are two cases to consider.

Case 1:  There is a column $\xx$ of $p$ that does not belong
to any child  of $p$ and is a linear combination of the rest of the
columns of $p$ (for instance this
happens when $p$ has at most one child).  In this case we
choose a vector $\tilde{\xx} \in Per(p) \setminus \mathrm{span}(p)$
which is guaranteed to exist by Lemma \ref{lem:enough}.  Then
adding an infinitesimal multiple of $\tilde{\xx}$ to $\xx$ increases
$\dim(p)$ without increasing the rank
of any of the matrices in $\mathcal{M}(p)$,
and hence does not change the fact that
$\mathbf{X}$ satisfies the minors of $P_{\Gamma}$.

Case 2:  Our element $p$ has a left child $p_1$ and a right child
$p_2$, but none of the columns of $p$ that are not in $p_1$ or
$p_2$ can be written as a linear combination of the rest of the
columns of $p$. We cannot add a vector $\tilde{\xx} \in Per(p)$ to
any part of $p$ which will increase $\dim(p)$ without risking the
increase of $\dim(p_1)$ or $\dim(p_2)$. We let $q$ be the common
child of $p_1$ and $p_2$, and if there is no such child we let $q
= \emptyset$. Now there exists a column $\xx$ of $p_1$ that is not
in $\mathrm{span}(q)$. This is clear when $q = \emptyset$, and
otherwise this follows from the minimality assumption on $p$. Now
we choose a vector $\tilde{\xx} \in Per(p) \setminus
\mathrm{span}(p)$ which is almost parallel to $\xx$, and we assume
that both vectors have the same norm. We let $B$ be a basis of
$\CC^m$ that contains the columns of $q$ (the columns of 
$q$ are linearly independent since $q$ is in at least the third
row of $\mathcal{Q}_\Gamma(\mathbf{X})$) as well as a basis for
$\mathrm{span}(p_1)$, and in particular $\xx$. This implies that
$\tilde{B}  = B \setminus \{\xx\} \cup \{ \tilde{\xx} \} $ is also
a basis for $\mathbb{C}^m$. We let $\mathbf{T}$ be the change of
basis matrix from $\tilde{B}$ to $B$. Now assuming that $p_2 =
[a,b]$, we perturb $\mathbf{X}$ and obtain
 
$$\mathbf{X}'= \left[
\begin{array}{c|c}
\mathbf{T} \cdot \mathbf{X}[1,a-1] & \mathbf{X}[a,n]
\end{array} \right].$$

\noindent Since $\tilde{\xx}$ is almost parallel to
$\xx$ and both vectors have the same norm, the linear
transformation $\mathbf{T}$ is small in the sense that it is close
to the identity matrix in the Euclidean topology. Furthermore,
this perturbation increases $\dim(p)$ by one, and any submatrix $q
\in \mathcal{M}(p)$ will not change its rank.  The rank of any 
submatrix $q
\geq p$ with $q \notin \mathcal{M}(p)$ increases by at
most one, and hence $\dim(q) \leq D(q)$ after perturbing by
$\mathbf{T}$.  Finally, a submatrix which does not contain $p$ is
either unchanged or is changed by applying an element of
$\gl_m(\mathbb{C})$ which does not alter the rank. This implies
that our new perturbed matrix is in $\mathcal{V}(P_\Gamma)$ and
completes the proof that we can always make perturbations to
improve the ranks of rank-deficient submatrices.
\end{proof}

\begin{thm} \label{thm:prime}  Let $K$ be a field of characteristic zero.  Then
$P_{\Gamma}$ is a prime ideal.
\end{thm}

\begin{proof}
First suppose that $K = \CC$.  Corollary \ref{cor:rad} says that
$P_\Gamma$ is radical and Propositions \ref{prop:irred} and
\ref{prop:perturb} imply that $\mathcal{V}(P_\Gamma)$ is
irreducible, hence $P_\Gamma$ is prime by the Nullstellensatz.
Now we apply the Lefschetz principle to deduce that $P_\Gamma$ is
prime over an arbitrary field $K$ of characteristic zero.  For
this suppose there are $f, g$ in $K[x_{ij}]$ with $fg \in
P_\Gamma$ but $f,g \notin P_\Gamma$. Then $fg \in P_\Gamma$ 
but $f,g \notin P_\Gamma$  over the field $\QQ(\{c_\alpha\})$ where
$\{ c_\alpha \}$ is the finite set of coefficients of $f$ and $g$.
Since $\CC$ has infinite transcendence degree over $\QQ$ and is
algebraically closed, and these fields have characteristic zero, 
$\QQ(\{c_\alpha\})$ can be embedded as a
subfield of $\CC$. The images of $f$ and $g$ under this 
embedding will show that $P_\Gamma$ is not prime over $\CC$.  
This is a contradiction.
\end{proof}

\section{Higher Dimensional Adjacent Minors}

Let $\mathbf{m} = (m_1, \ldots m_d) \in \ZZ^d$ with all $m_j \geq
2$ and  let $X_{\mathbf{m}}$ be the generic $d$-dimensional $m_1
\times \cdots \times m_d$ matrix with entries $x_{i_1, \ldots ,
i_d}$. Throughout this section we will call any integer vector
$\mathbf{u} = (u_1, \ldots, u_d)$ {\sl even} if $\sum u_j$ is
even, and ${\sl odd}$ otherwise.

\begin{defn}
Let $\mathbf{i} = (i_1, \ldots, i_d)$ be an integer vector with
$1 \leq i_j \leq m_j-1$ for all $j$.
A multidimensional adjacent 2-minor is a binomial of degree $2^{d-1}$
of the form

$$
\prod_{ \mathbf{\epsilon} \in \{0,1\}^d \atop
\mathbf{\epsilon} \,\, even} x_{\mathbf{i} + \mathbf{\epsilon}} -
\prod_{ \mathbf{\epsilon} \in
\{0,1\}^d \atop \mathbf{\epsilon} \,\, odd}
x_{\mathbf{i} + \mathbf{\epsilon}}.
$$

\noindent  Furthermore we let $I_{\mathbf{m}}(2)$ be the ideal in
$K[x_{\mathbf{i}}]$ generated by all the multidimensional adjacent
2-minors.
\end{defn}

The ideal $I_{\mathbf{m}}(2)$ generalizes the ideals $I_{mn}(2)$
of $2 \times 2$ adjacent minors from Section 2.  The set of
vectors $\{u-v \,:\, \xx^u - \xx^v \mbox{ is an adjacent 2-minor}
\}$ is a basis for the lattice of $d$-dimensional $m_1 \times
\cdots \times m_d$ integral matrices with all line sums equal to
zero \cite{HoSul}.  A \emph{line sum} of a matrix 
with entries $u_{\mathbf{i}}$ is any sum of the form
$$ \sum_{i_j =1}^{m_j} u_{\mathbf{i}}.
$$

\noindent This is actually only a very special case of the types
of marginals which one may compute of multidimensional matrices.
In fact, any marginal computation of a multidimensional matrix
leads naturally to a lattice of integer matrices with all
marginals equal to zero.  From this lattice, we can extract a
lattice basis of generalized adjacent minors \cite{HoSul}, and
construct an ideal of generalized adjacent minors. The general
results on lattice basis ideals in \cite{HS} imply that every
minimal prime of these ideals of adjacent minors is of the form in
equation (\ref{primeform}), so we only need to determine the
variables which appear in each minimal prime.

The similarity between the $2 \times 2$ adjacent minors for two
dimensional matrices and the higher dimensional adjacent minors we
describe in this section is somewhat misleading.  One important
difference is that the higher dimensional minors do not describe
rank conditions on tensors, so the linear algebra arguments which
we applied in Sections 3 and 4 no longer succeed. This problem
aside, one might still hope that the partition of variables which
arises in the description of the minimal primes of ideals of  higher
dimensional adjacent minors might still provide a decomposition of
the multidimensional matrix into rectangular chambers and their
boundaries. Unfortunately, this hope is far from the true
description of the minimal primes.  In this section, we describe
the minimal primes in a few special instances, showcasing the
increasing complexity which arises in higher dimensions.

\begin{example}\label{ex:223}
Let $d =3$ and $\mathbf{m} = (2,2,3)$.  The ideal of
multidimensional adjacent minors is

$$I_{2,2,3}(2)=  \langle \underline{x_{111}x_{122}x_{212}x_{221}}-x_{112}x_{121}x_{211}x_{222},
\underline{x_{112}x_{123}x_{213}x_{222}}-x_{113}x_{122}x_{212}x_{223} \rangle.
$$

\noindent  If we choose a term order which selects the underlined
terms as the leading terms, these leading monomials are relatively prime,
and hence this ideal is a
radical complete intersection.  The five minimal primes of
$I_{2,2,3}(2)$ are the ideals
$$ \langle x_{112},x_{122} \rangle, \langle x_{222},x_{122} \rangle,
\langle x_{222},x_{212} \rangle, \langle x_{212},x_{112} \rangle,
$$
\noindent and
$$
I_{2,2,3}(2): (\prod x_{ijk})^{\infty} = I_{2,2,3}(2) + \langle
x_{111}x_{123}x_{213}x_{221} - x_{113}x_{121}x_{211}x_{223}
\rangle.
$$
\end{example}

Generalizing Example \ref{ex:223} it is possible to give a
combinatorial description of the minimal primes of the ideal
of multidimensional adjacent minors whenever
$\mathbf{m} = (2, 2, \ldots, 2, m)$.
This is the content of the following theorem.

\begin{thm}\label{thm:22n}
Let $I_{\mathbf{m}}(2)$ be the ideal of adjacent 2-minors where
$\mathbf{m} = (2, 2, \ldots, 2, m)$. The minimal primes of
$I_{\mathbf{m}}(2)$ are of the form as in (\ref{primeform}) where
the set $S$ of variables is a collection of the pairs of variables
$x_{s_1, \ldots, s_{d-1}, j_i}$ and $x_{t_1, \ldots, t_{d-1},
j_i}$ chosen for each $j_i$ in the (possibly empty) set $J = \{2
\leq j_1 < \cdots < j_{\ell} \leq m-1 | j_i + 1 < j_{i+1} \}$ such
that the index vector of the first variable is even and the second
one is odd. Moreover, if we let $f_d(m)$ denote the number of
minimal primes of this ideal, then the function $f_d$ satisfies
the recurrence relation

$$f_d(m + 1) = f_d(m) + 4^{d-2} f_d(m-1),$$

\noindent with initial conditions $f_d(1) = f_d(2) = 1$.

\end{thm}
\begin{proof}
Note that if $P$ is a minimal prime of $I_{\mathbf{m}}(2)$ and
contains a variable $x_{s_1, \ldots, s_{d-1}, j_i}$ whose index
set is even (or odd) then it must contain some other variable
$x_{t_1, \ldots, t_{d-1}, j_i}$ whose index set is odd (or
respectively even). Moreover $P$ cannot contain another variable with
last index $j_i$ because this would contradict the minimality of
$P$. Any adjacent $2$-minor that contains these two variables must
contain two other variables of opposite parity with either a last
index $j_i + 1$ or $j_i - 1$, therefore no variable of this form
is needed in $P$. A similar reasoning implies that the variables
with last index $1$ or $m$ do not appear in $P$ either. This shows
that every minimal prime has the desired form. To see that every
ideal of the form we described is a minimal prime one needs merely
note that there are no containment relations between these ideals.

Now we prove the recurrence relation. Let $P$ be a minimal prime
arising from the sequence $J$. If the last index in $J$ is not
equal to $m-1$, then the sequence $J$ and the choice of variables
provides a minimal prime for $I_{\mathbf{m}}(2)$ where $\mathbf{m}
= (2, 2, \ldots, 2 ,m-1)$.  If the last index in $J$ is equal to
$m-1$, then removing it from the sequence (and the corresponding
variables from $S$) produces a minimal prime $Q$ for
$I_{\mathbf{m}}(2)$ where $\mathbf{m} = (2, 2, \ldots, 2 ,m-2)$.
There are precisely $4^{d-2}$ minimal primes $P$ that would give
rise to $Q$ since there are $4^{d-2}$ possible pairs of variables
with last index $m-1$ and having opposite parity.
\end{proof}

Aside from Theorem \ref{thm:22n} we do not know of any general
characterization of the minimal primes of these ideals of higher
dimensional adjacent minors.  We conclude this section with an
example which shows that these minimal primes do not have the same
appearance as in the two dimensional case, where the partition of
variables corresponded to rectangular subregions and their
boundaries.

\begin{example}
Let $\mathbf{m} = (3,3,3)$.  Then there are sixty-seven minimal
primes of $I_{\mathbf{m}}(2)$ which fall into nine symmetry
classes modulo the natural symmetry of the cube. In the table
below, we display the set of variables $S$ which appear in the
representative minimal primes, as well as the number of minimal
primes in a given symmetry class, and the degree of the
corresponding prime ideal.

$$ \begin{array}{|c|c|c|}
\hline S &  \mbox{size} &  \mbox{degree} \\
\hline \emptyset & 1 &  2457 \\
\hline \{x_{221}, x_{222}, x_{223} \} & 3 & 1 \\
\hline \{x_{121}, x_{122}, x_{123} \} & 12 &  81 \\
\hline \{x_{121}, x_{122}, x_{123}, x_{223}, x_{323} \} &  12  & 12 \\
\hline \{x_{121}, x_{122}, x_{123}, x_{232},x_{332} \} & 12 & 12 \\
\hline \{x_{122},x_{322},x_{211},x_{213},x_{231},x_{233} \} & 3 &1
\\
\hline \{x_{121}, x_{122}, x_{123}, x_{321}, x_{322}, x_{323} \}
& 6 & 1 \\
\hline \{x_{121}, x_{122}, x_{123}, x_{312},x_{322},x_{332}\} & 6
& 1 \\
\hline \{x_{121}, x_{123}, x_{232},x_{332}, x_{212},x_{312} \} &
12 & 1 \\
 \hline
\end{array}
$$
\end{example}

\bibliography{adjacent}
\bibliographystyle{plain}

\end{document}